\documentclass[pdflatex,sn-aps]{sn-jnl}

\usepackage{graphicx}%
\usepackage{multirow}%
\usepackage{amsmath,amssymb,amsfonts}%
\usepackage{amsthm}%
\usepackage{mathrsfs}%
\usepackage[title]{appendix}%
\usepackage{xcolor}%
\usepackage{textcomp}%
\usepackage{manyfoot}%
\usepackage{booktabs}%
\usepackage{algorithm}%
\usepackage{algorithmicx}%
\usepackage{algpseudocode}%
\usepackage{listings}%

\usepackage{natbib}

\theoremstyle{thmstyleone}%
\theoremstyle{thmstyletwo}%

\theoremstyle{thmstylethree}%

\begin{document}

\title{Globalizing manifold-based reduced models for equations and data}


\author[1]{\fnm{Bálint} \sur{Kaszás}}\email{bkaszas@ethz.ch}

\author[1]{\fnm{George} \sur{Haller}}\email{georgehaller@ethz.ch}

\affil*[1]{\orgdiv{Institute for Mechanical Systems}, \orgname{ETH Zürich}, \orgaddress{\street{Leonhardstrasse 21}, \city{Zurich}, \postcode{8092}, \country{Switzerland}}}
\maketitle
\begin{abstract}

{\textbf{Abstract}}:
 {One of the very few mathematically rigorous nonlinear model reduction methods is the restriction of a dynamical system to a low-dimensional, sufficiently smooth}, attracting invariant manifold. Such manifolds are usually found using local polynomial approximations and, hence, are limited by the unknown domains of convergence of their Taylor expansions. 
To address this limitation, we extend local expansions for invariant manifolds via Padé approximants, which re-express the Taylor expansions as rational functions for broader utility. This approach significantly expands the range of applicability of manifold-reduced models, enabling reduced modeling of global phenomena, such as large-scale oscillations and chaotic attractors of finite element models. We illustrate the power of globalized manifold-based model reduction on several equation-driven and data-driven examples from solid mechanics and fluid mechanics.
\end{abstract}

\section{Introduction}
Reduced-order modeling is a widespread technique that seeks to simplify high-dimensional nonlinear systems while retaining their essential dynamical features. Among reduced-order modeling procedures, manifold-based methods have been steadily gaining momentum. This can be largely attributed to the prevalence of data-driven approaches that successfully build on the fundamental manifold hypothesis \cite{baldi89,gorban2018, floryan2022data} of machine learning.  

Center manifold reduction \cite{guckenheimerNonlinearOscillationsDynamical1983}, geometric singular perturbation theory \cite{fenichelPersistenceSmoothnessInvariant1971, kuehnMultipleTimeScale2015} and inertial manifold theory \cite{foias88, koronaki2024} all rely on the existence of low-dimensional attracting invariant manifolds in the phase space of a dynamical system. These methods constitute mathematically rigorous examples of nonlinear model reduction and yield truly predictive models. However, for systems encountered in practice that are not close to bifurcations, invariant manifolds can only be realistically constructed when they emanate from a known robust stationary state, such as a hyperbolic fixed point. 

In those cases, seeking the invariant manifolds perturbatively and expressing them as local Taylor expansions at the known stationary state is justified. Traditionally, only stable and unstable manifolds, as continuations of the stable and unstable subspaces of the linearization, were approximated in this fashion \cite{guckenheimerNonlinearOscillationsDynamical1983}. The recent theory of spectral submanifolds (SSMs) extends this approach to arbitrary nonresonant spectral subspaces of the linearized system \cite{hallerNonlinearNormalModes2016b}. In particular, SSMs are now known to exist as smooth continuations of stable subspaces (like-mode SSMs) and of subspaces spanned by stable and unstable modes (mixed-mode SSMs) \cite{hallerNonlinearModelReduction2023a}. This model reduction approach has been used in a broad range of physical settings to deduce very low-dimensional, mathematically justified polynomial models \cite{hallerNonlinearNormalModes2016b, hallerNonlinearModelReduction2023a}.

SSM reduction has been successfully applied to obtain accurate reduced models of nonlinear vibrations observed in high-dimensional finite element models \cite{ jainHowComputeInvariant2022b, bettini2024} and experiments \cite{cenedeseDatadrivenModelingPrediction2022, cenedese2022mech}, multistable fluid flows \cite{kaszasDynamicsbasedMachineLearning2022,kaszasCapturingEdgeChaos2024}, chaotic systems \cite{liuDatadrivenModelingForecasting2024}, fluid-structure interaction problems \cite{xuDatadrivenModellingRegular2024, axasFastDatadrivenModel2023} and control of soft robots \cite{aloraPracticalDeploymentSpectral2023a}. In an equation-driven setting, SSM reduction starts with the solution of an invariance equation through local Taylor expansions \cite{cabreParameterizationMethodInvariant2003a,haroParameterizationMethodInvariant2016, jainHowComputeInvariant2022b}. Data-driven SSM reduction \cite{cenedeseDatadrivenModelingPrediction2022} also uses polynomial expansions to obtain an approximation for SSMs. 


SSMs are ideal tools for model reduction because their existence, uniqueness, and smoothness properties are precisely known. Specifically, if the governing equations are analytic, SSMs of attracting fixed points are guaranteed to also possess convergent Taylor series near the fixed point. However, the domain of convergence is generally unknown. 

This has been a general limitation of invariant manifold-based reduction methods, restricting the range of applicability of the resulting reduced-order model to an a priori unknown domain. Importantly, the convergence-limiting singularity of the Taylor expansion is not a physical singularity of the dynamical system, and hence, invariant manifolds continue to exist even outside the domain of convergence of the Taylor expansions for those manifolds. Therefore, the breakdown of convergence arises without any prior indication. This major limitation of local approximation methods of invariant manifolds restricts the user to potentially small physical domains. 

Here, we overcome this limitation by extending the local information contained in the Taylor series to considerably larger domains via the process generally known as analytic continuation \cite{rudinRealComplexAnalysis1987}. Among commonly used methods of analytic continuation, we focus on Padé approximants \cite{bakerPadeApproximants1996}. Padé approximants are rational functions whose Taylor expansion around a point coincides with that of the original function up to a given order but can represent the original function more efficiently and globally. 

Padé approximants have been used in theoretical physics and applied mathematics for summing divergent series \cite{benderAdvancedMathematicalMethods1999} with applications in cosmology \cite{alhoGlobalDynamicsInflationary2015a}, quantum electrodynamics \cite{zinn-justinConvergencePadeApproximants1974}, fluid dynamics \cite{razafindralandyNumericalDivergentSeries,andrianovPadeApproximantsTheir2021} and solving the Helmholtz equation \cite{bonizzoniLeastSquaresPadeApproximation2020}. Closest in spirit to our present work is the use of Padé approximants for center manifold reduction \cite{sinouMethodsReduceNonlinear2004, sinouNonlinearStabilityAnalysis2004}. The latter use is, however, restricted to equation-driven model reduction near Hopf bifurcations. 

In the context of SSM-reduced models, we must consider generalizations of Padé approximants to account for multivariate functions describing the parametrization and the reduced dynamics. Since two-dimensional non-resonant SSMs are typical in applications, we focus on the bivariate and univariate cases. 

 {In the following Section \ref{sec:results}, we present our results on using Padé approximants to construct reduced models on global SSMs (gSSMs). We discuss four examples, including the Kolmogorov flow \cite{chandlerInvariantRecurrentSolutions2013}, a nonlinear von Kármán beam in periodic and chaotic regimes \cite{jainHowComputeInvariant2022b, liuDatadrivenModelingForecasting2024}, and a data-driven model of an inverted flag experiment \cite{xuDatadrivenModellingRegular2024}. The mathematical details of SSMs and Padé approximants are discussed in Section \ref{sec:methods}. The Supplementary Material contains further applications and examples.}

\section{Results}\label{sec:results}
\subsection{Spectral submanifolds and Padé approximants}

We consider an $n-$dimensional nonlinear dynamical system
\begin{align}\label{eq:dynsys}
    \dot{\mathbf{x}} &= \mathbf{A}\mathbf{x} + \mathbf{f}(\mathbf{x}) + \varepsilon \mathbf{f}_{\text{ext}}(\mathbf{x}, t), \\ \mathbf{x} &\in \mathbb{R}^n, \quad \mathbf{A} \in \mathbb{R}^{n \times n}, \qquad \mathbf{f} \in \mathcal{O}(|\mathbf{x}|^2), \nonumber \\
    n &\geq 1, \qquad 0\leq \varepsilon \ll 1 \nonumber
    \end{align}
where $\mathbf{f}_{\text{ext}}(\mathbf{x},t) = \mathbf{f}_{\text{ext}}(\mathbf{x},t + T)$ represents time-periodic external forcing. We assume that the nonlinearity $\mathbf{f}$ and the forcing $\mathbf{f}_{\text{ext}}(\mathbf{x},t)$ are also analytic functions of $\mathbf{x}$. We assume that, for $\varepsilon=0$, the fixed point $\mathbf{x}=0$ is hyperbolic and the spectrum of $\mathbf{A}$ is nonresonant. The slow spectral subspace is denoted $E$, and is defined as the $d$-dimensional real subspace spanned by the eigenvectors of $\mathbf{A}$ associated to its $d$ eigenvalues closest to zero and hence is an attracting, slow invariant manifold for the linearized dynamics. 

We focus here on spectral submanifolds (SSMs), which are the nonlinear continuations of spectral subspaces under the
addition of the nonlinear terms in \eqref{eq:dynsys}. Although there are multiple invariant manifolds tangent to the spectral subspace $E$, there is a unique smoothest one, called the primary spectral submanifold \cite{hallerNonlinearNormalModes2016b} denoted $\mathcal{W}(E)$. If the fixed point $\mathbf{x}=0$ is stable, then $\mathcal{W}(E)$ is even known to be analytic. Due to the slowness of the subspace $E$, the SSM, $\mathcal{W}(E)$, is an attracting slow manifold for the autonomous system \eqref{eq:dynsys}. The members of the invariant manifold family with reduced smoothness are called fractional (or secondary) SSMs \cite{hallerNonlinearModelReduction2023a}.
For completeness, the necessary assumptions for the existence and uniqueness of primary SSMs \cite{hallerNonlinearNormalModes2016b, hallerNonlinearModelReduction2023a, hallerNonlinearModelReduction2024} are recalled in Section \ref{sec:ssms}.

 {In the autonomous limit with $\varepsilon=0$, the $d-$dimensional ($d\leq n$) primary SSM, $\mathcal{W}(E)$, can locally be represented as the image of a parametrization $\mathbf{W}:U\subset \mathbb{R}^d \to \mathbb{R}^n$, over some open set $U\subset \mathbb{R}^d$ as
\begin{equation}
\label{eq:param}
    \mathcal{W}(E) = \left\{\mathbf{x}=\mathbf{W}(\mathbf{p})\ | \   \mathbf{p} \in U \right\}\subset \mathbb{R}^n.
\end{equation} 
The reduced dynamics $\dot{\mathbf{p}} = \mathbf{R}(\mathbf{p})$, with $\mathbf{R}:U\to \mathbb{R}^d$ are conjugate to \eqref{eq:dynsys}, i.e., $\mathcal{W}(E)$ is invariant under the time evolution of \eqref{eq:dynsys} and therefore satisfies the invariance equation 
\begin{align} \centering
    \label{eq:inveq_aut}
     \mathbf{A}\mathbf{W}(\mathbf{p}) + \mathbf{f}(\mathbf{W}(\mathbf{p})) =
    D\mathbf{W}(\mathbf{p})\dot{\mathbf{p}}.
\end{align}
We refer to Section \ref{sec:ssms} for a discussion on SSMs of the nonautonomous system with $\varepsilon>0$.  }

 {We solve Eq. \eqref{eq:inveq_aut} by representing the parametrization of $\mathcal{W}(E)$ and its reduced dynamics as a power series truncated to some order $N$, i.e., 
\begin{align}
    \label{eq:w_r_multiindex}
    \mathbf{W}^N(\mathbf{p}) &=  \sum_{|\mathbf{k}|=0}^N\mathbf{W}_{\mathbf{k}}\mathbf{p}^\mathbf{k},\\
    \mathbf{R}^N(\mathbf{p}) &= \sum_{|\mathbf{k}|=0}^N\mathbf{R}_{\mathbf{k}}\mathbf{p}^\mathbf{k}. \nonumber
\end{align}}

 {We define the multi index $\mathbf{k}=(k_1, ..., k_d)$ and $|\mathbf{k}|=k_1 + k_2 + ... + k_d$, so that $\mathbf{p}^\mathbf{k} = p_1^{k_1}p_2^{k_2}...p_d^{k_d}$ refers to a scalar monomial of the components of $\mathbf{p}$ with total order $|\mathbf{k}|$. The coefficients $\mathbf{W}_\mathbf{k}$ and $\mathbf{R}_\mathbf{k}$ are vectors in $\mathbb{R}^n$ and $\mathbb{R}^d$, respectively, for all $\mathbf{k}$. } 

 {The coefficients $\mathbf{R}_\mathbf{k}$ depend on the style of parametrization used. In the graph style parametrization, the reduced coordinates are obtained as projections onto the spectral subspace $E$, while in the normal form style parametrization, non-resonant terms are set to zero. The difference between these two choices is explained in more detail by, e.g., \cite{jainHowComputeInvariant2022b,haroParameterizationMethodInvariant2016, stoychevFailingParametrizationsWhat2023}. }

Since the primary autonomous SSM is analytic, there is a domain of convergence denoted as $U_0\subset U \subset \mathbb{C}^d$, for which the $N\to \infty$ limits exists, i.e.,
\begin{equation}
    \lim_{N\to\infty} \mathbf{W}^N(\mathbf{p}) = \mathbf{W}(\mathbf{p}), \quad \forall \mathbf{p} \in U_0.
\end{equation}
For a system whose slowest mode is oscillatory and not in resonance with higher modes, the optimal model reduction is achieved by a two-dimensional SSM tangent to a single oscillatory eigenspace of the autonomous problem. In that case, we can select $\mathbf{p}= (p, \bar{p})^T$ with $p\in \mathbb{C}$.  {With the normal form style parametrization, the reduced dynamics only contain near-resonant terms of the form $p^{k+1}\bar{p}^k$ and $p^k\bar{p}^{k+1}$, and it is conveniently expressed in polar coordinates \cite{jainHowComputeInvariant2022b}. We introduce}
\begin{equation}
    p = \rho e^{i \theta}, \quad \bar{p} = \rho e^{-i \theta},
\end{equation}
which allows us to write the SSM-reduced dynamics as
\begin{align}\label{eq:rho_theta_aut}
    \dot{\rho} &= \kappa(\rho)\rho, \\
    \dot{\theta} &= \omega(\rho).  \nonumber
\end{align}
The functions $\kappa(\rho)$ and $\omega(\rho)$ represent the instantaneous dependence of the damping rate on the amplitude and the frequency on the amplitude, respectively. These functions are obtained from the coefficients of the autonomous reduced dynamics, $\mathbf{R}_\mathbf{k}$. Their Taylor expansions are of the form
\begin{align}\label{eq:omegakappa}
    \kappa(\rho) = \sum_{n=0}^\infty \kappa_n \rho^{2n}, \quad
    \omega(\rho) = \sum_{n=0}^\infty \omega_n \rho^{2n}.
\end{align}

The software package \texttt{SSMTool} \cite{jainHowComputeInvariant2022b, ssmtool_software} returns the Taylor coefficients $\mathbf{W}_{\mathbf{k}}$, $\kappa_k$ and $\omega_k$ up to arbitrary orders. However, the expansions \eqref{eq:omegakappa} only converge as long as the amplitude $\rho$ corresponding to the physical response of the system is inside the domain of convergence. 

For the expansions in \eqref{eq:omegakappa}, the domain of convergence is the interior of a disk of radius $R$, the radius of convergence. As a corollary of the analyticity of holomorphic functions \cite{rudinRealComplexAnalysis1987}, that disk contains the nearest singularity of the complex extension of the function. This result also generalizes to the multivariate setting, although defining the domain of convergence is more complicated \cite{scheidemannIntroductionComplexAnalysis2023a}. In addition, in contrast to the scalar case, singularities of multivariate functions are never isolated. 

The a priori unknown domain of convergence of their local Taylor expansions represents a limitation of SSM-reduction approaches. Although solutions of the invariance equation \eqref{eq:inveq_periodic} for SSMs exist up to any order $N$, their formal sum $\mathbf{W}^N$ has little to do with the primary SSM outside the domain of convergence $U_0$, even if $\mathbf{W}$ is well-defined on $U\backslash U_0$. 

This fundamental limitation has forced most invariant manifold studies to focus on deriving reduced-order models under small perturbations near a fixed point. This, however, impedes predicting the system's response to large perturbations or the discovery of steady states far away from the known fixed points.

However, SSMs, as invariant manifolds, are known to extend globally in the phase space, as long as the flow map of the dynamical system \eqref{eq:dynsys} remains well defined on them for all times. A straightforward extension of the locally known parametrization is to evolve a set of initial conditions from inside the domain of convergence globally under the flow map. This, however, assumes that the flow map of the full high-dimensional system is known in detail, making the reduced-order model redundant. This envisioned globalization is only partly achieved by the data-driven construction of SSMs \cite{cenedeseDatadrivenModelingPrediction2022}, which starts from a limited number of trajectories of \eqref{eq:dynsys} and finds the observed invariant manifold using regression. 

In this work, we propose a different approach to extend the range of applicability of SSM-reduced models. As long as an analytic function is known on some open domain, fundamental results in complex analysis guarantee that this knowledge can be extended to the entire domain of analyticity of the function, possibly using a different representation of the function. Switching to such a representation, other than the Taylor series \eqref{eq:w_r_multiindex}, is known as analytic continuation \cite{rudinRealComplexAnalysis1987}, a powerful technique to {\em globalize} the maps $\mathbf{W}(\mathbf{p})$ and $\mathbf{R}(\mathbf{p})$.

\subsection{Globalization of invariant manifolds via Padé approximation}
A well-known method for analytic continuation is the Padé approximation \cite{bakerPadeApproximants1996}, which has had numerous applications in theoretical physics and engineering. To describe this procedure, we introduce a multivariate rational function of the form,  {using the same notation as in Eq. \eqref{eq:w_r_multiindex}},
\begin{align}
    [N/M](\mathbf{z}) = \frac{\sum_{|\mathbf{k}|=0}^N a_{\mathbf{k}}\mathbf{z}^\mathbf{k}}{\sum_{|\mathbf{k}|=0}^M b_{\mathbf{k}}\mathbf{z}^\mathbf{k}}, \quad b_{\mathbf{0}} = 1, \quad \mathbf{z}\in \mathbb{R}^\ell,
    \label{eq:2dpade}
\end{align}
where the orders of the numerator and denominator are fixed constants $N,M$. This formulation covers the cases of univariate ($\ell=1$) and multivariate  {$(\ell\geq2)$} functions as well. 

We represent solutions of the invariance equation \eqref{eq:inveq_periodic} as rational functions of the form  \eqref{eq:2dpade}, which is achieved by requiring that the Taylor series of \eqref{eq:2dpade} matches that of the function to be approximated. Rational functions of the form \eqref{eq:2dpade} have the advantage that they can incorporate singularities that would otherwise limit the convergence of Taylor series. 
Based on these observations, we seek to extend SSM-reduced models using Padé approximants. We call the extended representations of an SSM the {\em global SSM} (gSSM).

Increasing the orders of the numerator and the denominator in \eqref{eq:2dpade}, one expects that the approximants provide better approximations for the functions $\mathbf{W}(\mathbf{p})$, $\mathbf{R}(\mathbf{p})$, $\kappa(\rho)$, and $\omega(\rho)$. In most practical applications, diagonal approximants (i. e., those with numerators and denominators of the same order) have proven to be the most effective.  {Moreover, diagonal Padé approximants of a univariate function are related to the continued-fraction representation of the function. }

For meromorphic functions with an a priori unknown number of poles, strict convergence is only guaranteed in measure by the theorems of \cite{nuttallConvergencePadeApproximants1970} and \cite{pommerenkePadeApproximantsConvergence1973}. These theorems cover the case of diagonal approximants and state that the sequence $[M/M](z)$ converges to the given function as $M\to \infty$ on bounded compact subsets of $\mathbb{C}$, except for sets of measure zero. Similar theorems also exist for the multivariate case \cite{goncarLocalConditionSinglavaluedness1974}. These exceptional sets correspond to zero-sets of the denominator in \eqref{eq:2dpade}. 

\subsection{Data-driven global reduced-order models}
If the equations of motion of a dynamical system are known, we will rely on the established theory of Padé approximants and invariant manifolds and use highly optimized computational methods to solve the invariance equation \eqref{eq:inveq_periodic} and construct the gSSM-reduced models. 

In many practical applications, however, the governing equations are only partially known or completely unknown, and yet a predictive reduced-order model is required. In such cases, we must rely directly on data-driven methods to approximate the gSSMs. In particular, the \texttt{SSMLearn} algorithm of \cite{cenedeseDatadrivenModelingPrediction2022} works on observations of trajectories to approximate the slow SSM $\mathcal{W}(E)$ locally. 

A common use case is when a single scalar observation $y(t)\in\mathbb{R}$ is recorded. In that case, by the Takens embedding theorem \cite{takensReconstructionTheoryNonlinear2010}, a $d-$dimensional attracting SSM can be embedded generically in the space of delayed measurements,
\begin{equation}
    \mathbf{y}(t) = (y(t), y(t-\tau), ..., y(t-(p-1)\tau))^T,
\end{equation}
for some time-lag $\tau$ as long as the number of delays $p$ is more than twice the dimension of the SSM $d$. The slow SSM can then be parametrized as a graph over its tangent space at the fixed point from which it emanates. The reduced coordinates $\boldsymbol{\eta}$ on the SSM are obtained by projecting the delay embedded measurements onto that $d-$dimensional tangent space, i.e., letting 
\begin{equation}
    \boldsymbol{\eta} = \mathbf{V}^T \mathbf{y},
\end{equation} 
where the columns of $\mathbf{V}$ span the tangent space. The tangent space can be efficiently approximated using the leading principal components of the delay-embedded measurements after discarding initial transients \cite{axasFastDatadrivenModel2023,xuDatadrivenModellingRegular2024}.

The parametrization of the SSM is obtained by regression using the observed trajectories and their reduced coordinates. The most straightforward choice is a polynomial basis, which is justified by the existence of a locally convergent Taylor expansion of the SSM and by the universal approximation property of polynomials \cite{rudinRealComplexAnalysis1987}. 

Motivated by the success of Padé approximants in enhancing the convergence properties of equation-driven models, we generalize the regression step of the data-driven SSM-reduction to rational approximants. Related approaches are rational interpolation and multi-point Padé approximation. While the former requires the approximant to fit the data exactly, the latter matches the asymptotic behavior of the function at multiple expansion points. 

Rational function regression \cite{austinPracticalAlgorithmsMultivariate2021} seeks a rational approximant of the form \eqref{eq:2dpade}. In addition, we enforce that all components of the vector function share the same denominator and that the denominator is never zero  {on the training data. Having a common denominator for all components of the vector function makes it simpler to avoid spurious singularities.} We elaborate on the steps of this regression task in Section \ref{sec:regression}.

\subsection{Example 1: Connecting orbit in Kolmogorov flow}
Our first example is the 2D Kolmogorov flow, governed by the Navier-Stokes equations in a periodic domain subject to spatially periodic forcing. The flow domain is given by $x,y \in [0,2\pi]$ and the time evolution of the vorticity $\omega(x,y)= \nabla \times \mathbf{u}$ is governed by the non-dimensionalized equation
\begin{align}
    \label{eq:kolmogorov}
    \frac{\partial \omega }{\partial t} = - (\mathbf{u}\cdot \nabla )\omega + \frac{1}{\text{Re}}\Delta \omega - L\cos Ly,
\end{align}
 {where $L$ denotes the forcing wave number. This influences the size of the large-scale flow structures, the bifurcations observed in \eqref{eq:kolmogorov} and the properties of the turbulent dynamics at high \text{Re} \cite{platt91}}. The laminar solution is a fixed point of \eqref{eq:kolmogorov} for all Reynolds numbers that can be written as 
\begin{equation}
    \omega_0(x,y) = -\frac{\text{Re}}{L }\cos Ly.
\end{equation}

Following \cite{chandlerInvariantRecurrentSolutions2013}, we set $L = 4$ and discretize the system using 576=24$\times$24 Fourier modes. The numerical implementation is based on \cite{farazmandTensorbasedFlowReconstruction2023} and results in a system of ODEs in the form \eqref{eq:dynsys} for the Fourier amplitudes $\hat{\omega}(k_x, k_y)$ with $n=576$. Further details on the implementation can be found in the Supplementary Material.

The laminar flow $\omega_0$ is already unstable for low Reynolds numbers, losing stability in a bifurcation around $\text{Re} = 9.1$. The state $\omega_0$ is connected to two new stable fixed points $\omega_{1,2}$ by its 1D unstable manifold. This unstable manifold coincides with the 1D slow SSMs of $\omega_{1,2}$, forming two heteroclinic orbits.

We fix $\text{Re}=11$ and consider the stable fixed point $\omega_{1}$. We compute the parametrized slow SSM of $\omega_1$ as i.e., $\hat{\omega} = \mathbf{W}(\xi)$. This computation is carried out automatically by \texttt{SSMTool}, which returns the coefficients of the Taylor expansion of $\mathbf{W}(\xi)$ and the reduced dynamics on the SSM as
\begin{align}
    \hat{\omega}(k_x, k_y) &= \mathbf{W}(\xi) = \sum_{k}\mathbf{W}_k \xi^k, \\
    \dot{\xi} &= R(\xi) = \sum_{k}R_k \xi^k.
\end{align}

\begin{figure}[H]
    \centering
    \includegraphics[width = 0.99\linewidth]{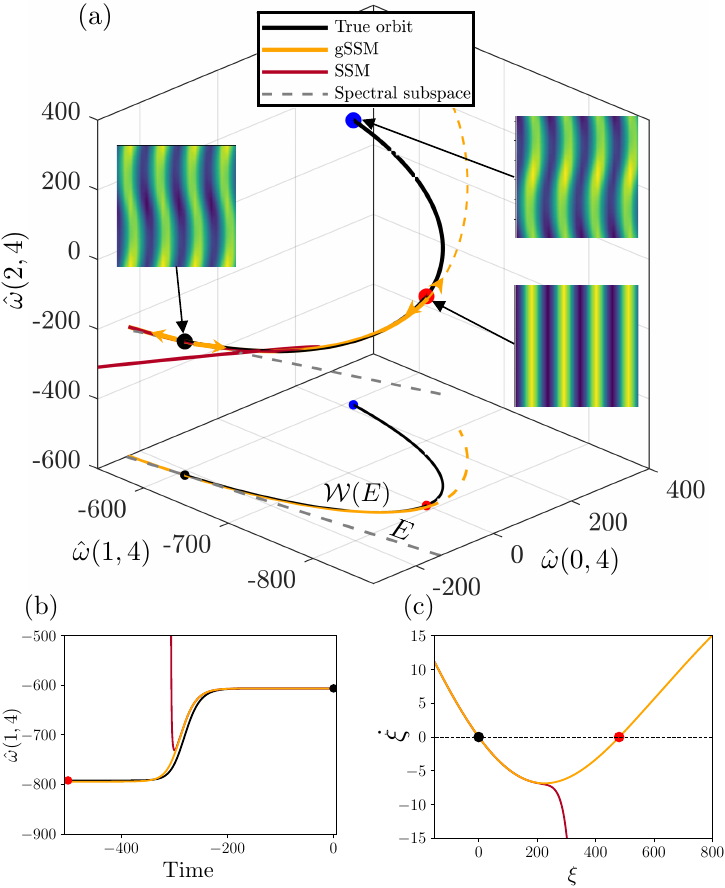}
    \caption{Heteroclinic orbits (black) of the Kolmogorov flow connecting $\omega_{1,2}$ and $\omega_{0}$. (a): Projection of the phase space onto  {three dominant Fourier modes} $(1,4)$, ($0,4$) and ($2,4$).  {We show the slow SSM $\mathcal{W}(E)$ (black), which is tangent to the spectral subspace $E$ (grey), its order-16 Taylor expansion (red) near the fixed point $\omega_{1}$, and the order [5/5] gSSM approximation (orange). The curves are also projected to the horizontal axes. The vorticity fields corresponding to the three fixed points are shown in the insets.  (b): A trajectory on the heteroclinic orbit obtained by backward integration and its SSM-reduced and gSSM-reduced counterparts. (c): SSM-reduced and gSSM-reduced dynamics. }}\label{fig:kolmogorov}
\end{figure}

We visualize the connecting orbit, i.e., the 1D slow SSM in a 3D slice of the 576-D phase space in Fig. \ref{fig:kolmogorov}a. The stable fixed points $\omega_{1,2}$ are marked with black and blue dots, and the unstable fixed point $\omega_0$ is red. The heteroclinic orbit, which is obtained by direct numerical integration of \eqref{eq:kolmogorov} (black curve), is approximated by the Taylor expansion of the SSM up to order 16. The domain of convergence is clearly limited, and it does not contain the unstable fixed point.

To construct a globalized slow SSM, we compute the Padé approximant of each component (Fourier mode) of the vector $\mathbf{W}(\xi)$. Although, in principle, spurious poles could arise for each component separately, we find that the diagonal approximants are well-behaved. An alternative approach would be to construct the vector Padé approximants \cite{sidiVectorExtrapolationMethods2017}, that is, to find a common denominator for all components of $\mathbf{W}(\xi)$ \cite{laneSteadyBoussinesqConvection2024}.

The componentwise computed $[5/5](\xi)$ Padé approximant approximates the heteroclinic orbit connecting $\omega_1$ and $\omega_0$ remarkably well, as shown by the orange curves in Fig. \ref{fig:kolmogorov}a.  {We also observe that the manifold $\mathcal{W}(E)$ can only be represented as a graph over $E$ for this segment, since the derivative $\frac{\partial}{\partial \xi} \mathbf{W}(\xi)$ diverges at a fold point near $\omega_0$. The parametrization, therefore, cannot be continued to capture $\omega_2$. Note, however, that the Taylor approximation diverges well before encountering this unremovable singularity of the graph-style parametrization.}

 {To verify the validity of the reduced-order models, predictions should be compared to trajectories of the full system. However, since the fixed point is stable, a nearby initial condition will leave its neighborhood along the heteroclinic orbit only in backward time. Therefore, we integrate the initial condition $\xi(0)=10^{-5}$, close to the stable fixed point, backward under the SSM-reduced and the gSSM-reduced dynamics.} 

Fig. \ref{fig:kolmogorov}b shows the SSM-reduced trajectory, which exhibits finite-time blowup and is only a reliable model near the stable fixed point. In contrast, the gSSM-reduced trajectory converges to $\omega_0$ in backward time and hence captures the global dynamics accurately.  

Fig. \ref{fig:kolmogorov}c shows the reduced dynamics $\dot{\xi}$ on the SSM and the gSSM. In addition to the trivial fixed point at $\xi=0$, the gSSM-reduced model contains the unstable fixed point $\omega_0$, given by the intersection with $\dot{\xi}=0$ and hence provides a robust reduced model of the system. In contrast, based on the Taylor approximated SSM, the existence of an unstable fixed point cannot be concluded. 

\subsection{Example 2: Von Kármán beam}\label{sec:vonkarman}
\begin{figure*}[h!]
    \centering
    \includegraphics[width = 0.99\linewidth]{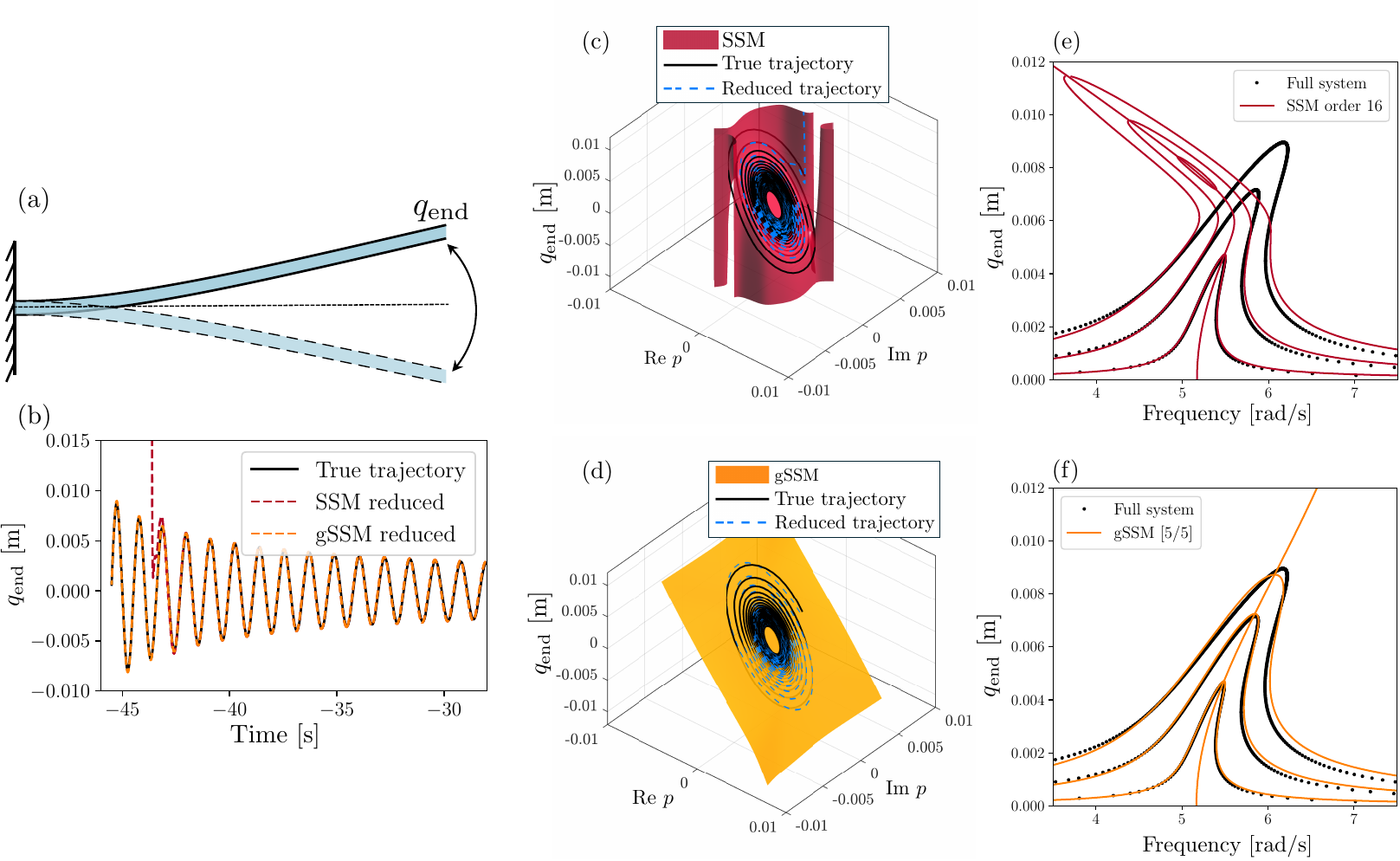}
    \caption{(a): von Kármán beam with clamped-free boundary conditions. (b): Trajectory of the unforced system (black) with its order-16 SSM approximation (red) and gSSM-approximation (orange) with a [5/5] and [5/4] Padé approximant. (c)-(d): Representation of the end point displacement $q_{\text{end}}$ and the full-order trajectory shown in (b) with the SSM and gSSM-reduced trajectories. (e)-(f): Forced response defined as the maximal end point displacement 
 due to a forcing amplitude $\varepsilon = 0.5, 1.7, 2.8$. Panel (e) shows the SSM-prediction and (f) shows the gSSM prediction. A supplementary animation showing the gSSM-prediction is available at \url{https://polybox.ethz.ch/index.php/s/ePwsPfDgHnlzlJ2}.}\label{fig:vonkarman_forced}
\end{figure*}

We now consider the model of a nonlinear von Kármán beam with clamped-free boundary conditions \cite{jainExactNonlinearModel2018}, shown in Fig. \ref{fig:vonkarman_forced}a. The beam has length $L=1 \ \text{m}$, height $h=1 \ \text{mm}$, width $b=0.1 \ \text{m}$, Young's modulus $E=70 \ \text{GPA}$, viscous damping rate $\alpha = 10^7 \ \text{Pa}\cdot\text{s}$ and density $2700 \text{kg/m}^3$. We use a finite element discretization with 10 elements, resulting in $30$ degrees of freedom, with a phase space of dimension $n=60$. The discretization is implemented in the repository \texttt{SSMTool} \cite{ssmtool_software}. 

We construct the slowest SSM of the undeformed configuration, which is tangent to the spectral subspace corresponding to the eigenvalues $\lambda_{1,2} = -0.0019 \pm 5.1681 i$. Due to the Taylor approximation, the autonomous SSM can only capture the decaying oscillations of the beam up to amplitudes of around $5$ mm, which is verified by backward integration of the SSM-reduced and the full models. This can be seen in Fig. \ref{fig:vonkarman_forced}b-c, which shows the autonomous trajectory exiting the domain of convergence of the Taylor series. This is evident for the high-order Taylor approximant, which exhibits finite-time blowup. 

In contrast, the gSSM model, globalized using a [5/5] Padé approximant for the reduced dynamics and a [5/4] approximant for the parametrization, remains well-behaved for even larger amplitudes. Due to spurious singularities in the parametrization, we chose the [5/4] approximant instead of [5/5]. Since the convergence of Padé-approximants is only guaranteed in measure, singularities coinciding with zero-sets of the denominator of \eqref{eq:2dpade} must be actively avoided before deploying the reduced-order model, as we demonstrate in the Supplementary Material. 

The backbone curve, i.e., the instantaneous relationship between the normal form amplitude $\rho$ and the frequency is given by $\omega(\rho)$ in \eqref{eq:rho_theta_aut}. The Taylor expansion of $\omega(\rho)$ and its $[5/5]$ Padé approximant are given as
\begin{align}
    \label{eq:vonkarman_omega}
    \omega(\rho) &= 5.16 +9.3\cdot 10^4 \rho^2 - 4.5\cdot 10^9 \rho^4 \nonumber \\
    &+ 4.2\cdot 10^{14} \rho^6 + O(\rho^8), \\
    [5/5](\rho) &= \frac{5.16 + 10^6 \rho^2 + 3.7\cdot 10^{10}\rho^4}{1 + 1.8\cdot 10^5 \rho^2 + 4.7\cdot \rho^4}.\nonumber 
\end{align}
 The Taylor coefficients in \eqref{eq:vonkarman_omega} are growing rapidly, and the alternating sign pattern of the coefficients suggests the convergence limiting singularity is along the imaginary axis. This can be inferred using the method of \cite{mercerCentreManifoldDescription1990b, dombSusceptibilityFerromagneticCurie1957}, which we specialize to our examples in the Supplementary Material.

When an external forcing $\varepsilon \mathbf{f}_{ext}\cos (\Omega t)$ is applied to the beam with a frequency $\Omega$ almost in resonance with the slowest eigenfrequency $\text{Im } (\lambda_{1,2})$, we can use the unforced SSM to make predictions about the forced response. To leading order, the reduced dynamics of the forced system are simply a perturbed version of those of the unforced one (see Section \ref{sec:ssms}). Furthermore, the forced response corresponds to the periodic orbits of the reduced model. For 2D SSMs, these are directly given by the equation \eqref{eq:forced_resp_eq} in Section \ref{sec:ssms}.

We compare the SSM-reduced forced response to the forced response of the full system, obtained by direct numerical continuation using COCO \cite{dankowiczRecipesContinuation2013}. Predictions with a high-order Taylor approximation of the SSM are shown in Fig. \ref{fig:vonkarman_forced}e. They accurately capture the forced response for small forcing amplitudes up to $\varepsilon = 0.5$, as initially reported by \cite{jainHowComputeInvariant2022b}. For higher amplitudes, the Taylor approximated reduced trajectory reaches its boundary of convergence, and the model breaks down. In contrast, the gSSM prediction using the [5/5] Padé approximant remains accurate for larger amplitudes as well. Note that the small errors in the peaks of the forced response curve are due to our initial assumption of a small forcing amplitude. Indeed, \eqref{eq:forced_resp_eq} technically holds only for small $\varepsilon$. 

\subsection{Example 3: Chaotic von Kármán beam}\label{sec:buckled_beam}

\begin{figure*}[h!]
    \centering
    \includegraphics[width = 0.99\linewidth]{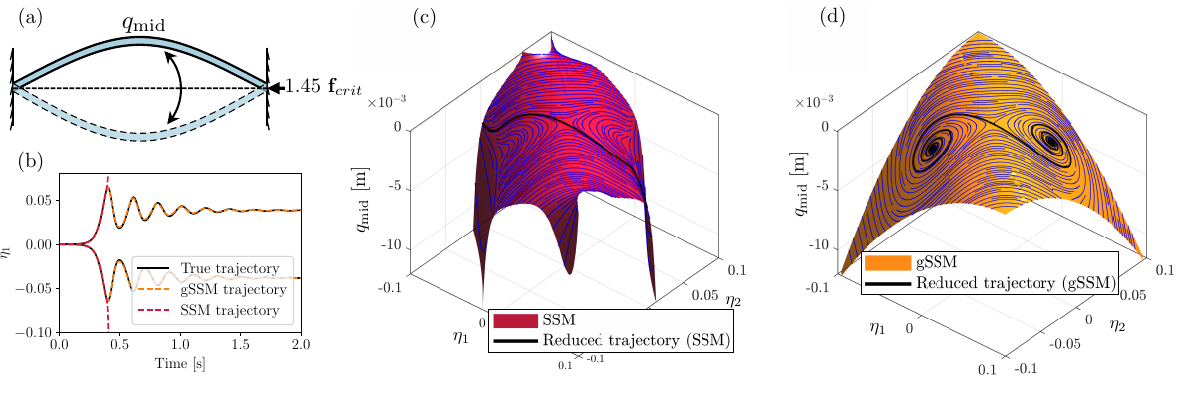}
    \caption{(a): Sketch of the beam in the buckled configuration, with no external forcing. (b): Trajectories in the unstable manifold of the unstable fixed point (black). Their order-18 SSM-reduced (red) and order [6/6] gSSM-reduced (orange) approximations are also shown. (c)-(d): The SSM and gSSM in the physical space with the direction field of the reduced dynamics indicated on the surface of the manifold. The predicted trajectories connecting the unstable fixed point to the stable fixed points are shown in black. }\label{fig:vonkarman_buckled}
\end{figure*}

In our next example, we construct a 2D mixed-mode SSM to characterize the chaotic behavior of a periodically forced buckled von Kármán beam, shown in Fig. \ref{fig:vonkarman_buckled}. We adopt the same finite element code used in Section \ref{sec:vonkarman}, developed by \cite{jainHowComputeInvariant2022b}, with pinned-pinned boundary conditions. To induce buckling, a compressive force is applied at the rightmost element, equal to $145\%$-of the critical value $f_{crit}=1.5 \ $ kN. This gives rise to a pair of stable fixed points. The same system, much closer to the critical point, was also analyzed using a data-driven model by \cite{liuDatadrivenModelingForecasting2024}. 

The spectral subspace of the linear part associated with the buckling instability is two-dimensional and has the corresponding eigenvalues $\lambda_{1} = 23.48$ and $\lambda_2 = -23.48$. The other modes are all stable and correspond to oscillatory dynamics. This eigenvalue configuration indicates a considerable departure from criticality and center manifold theory.

We use \texttt{SSMTool} to find an order-18 approximation of the 2D SSM of the unforced beam and its reduced dynamics. Since the SSM is tangent to a spectral subspace with real eigenvalues, we denote the reduced coordinates as $\boldsymbol{\eta}=(\eta_1, \eta_2)\in \mathbb{R}^2$. The reduced dynamics are given by 
\begin{align}
    \dot{\eta}_1 &= 23.48 \eta_1 -2800\eta_1^3 -1760\eta_2^3-7297\eta_1^2\eta_2 \nonumber \\
    &-6268 \eta_1\eta_2^2 + O(\boldsymbol{|\eta|}^4), \nonumber \\
    \dot{\eta}_2 &= -23.52 \eta_2 +2790 \eta_1^3 +1760\eta_2^3+7297\eta_1^2\eta_2 \\
    &+6268 \eta_1\eta_2^2 + O(\boldsymbol{|\eta|}^4). \nonumber 
\end{align}

Analysis of the reduced dynamics shows, however, that the fixed points born due to the buckling instability lie outside the domain of convergence. We show in Fig. \ref{fig:vonkarman_buckled}b the time series of a trajectory initialized on the unstable manifold of the fixed point. The SSM-reduced trajectory based on Taylor expansion blows up once the reduced trajectory exits the domain of convergence. In contrast, a comparable, [6/6] Padé approximant and the gSSM reduced-trajectory correctly captures the convergence to both of the buckled states.

This is even more apparent in Fig. \ref{fig:vonkarman_buckled}c-d, showing the image of the parametrization and the direction field of the reduced dynamics. The SSM-reduced model in Fig. \ref{fig:vonkarman_buckled}b predicts diverging, unphysical displacements for the mid-point of the beam. The gSSM-model in Fig. \ref{fig:vonkarman_buckled}c correctly identifies all fixed points and the orbits connecting them. We have, therefore, extended the SSM-reduced model obtained from local information around the unstable fixed point to a globally valid one. 

\begin{figure*}[h!]
    \centering
    \includegraphics[width = .99\linewidth]{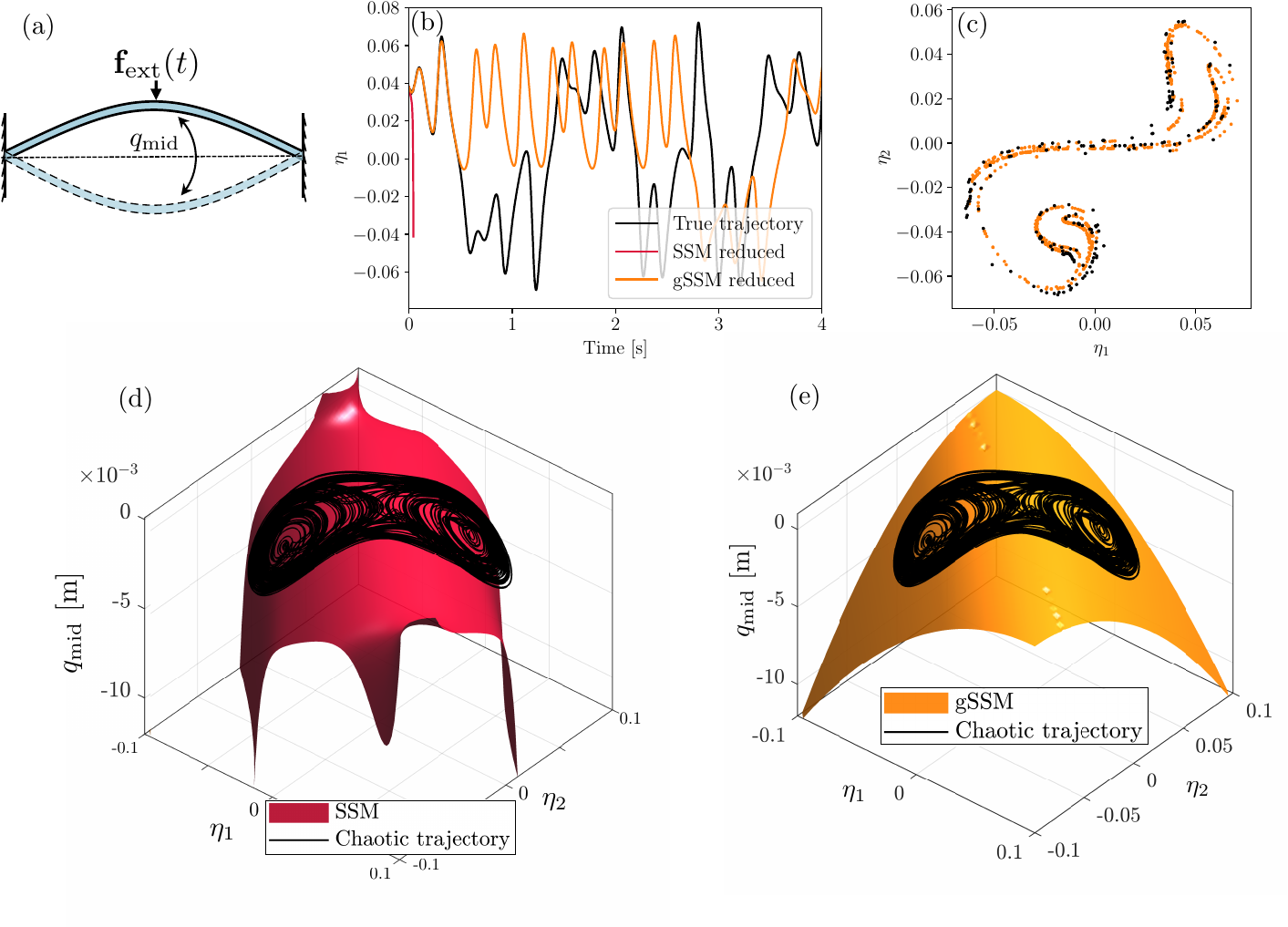}
    \caption{(a): Buckled von Kármán beam with periodic external forcing. (b): Time series of the reduced coordinate $\eta_1$ on a chaotic trajectory of the full system (black). Also shown are the SSM-reduced forced model, which diverges immediately (red), and the gSSM-reduced trajectory (orange). (c): Sampling of the Poincaré map of the gSSM-reduced model (orange) and the true system (black). (d)-(e): The autonomous SSM and gSSM with the chaotic trajectory of the full model. A supplementary animation comparing the SSM- and gSSM-predictions is available at \url{https://polybox.ethz.ch/index.php/s/ozX0r0Gx9X2Ryk2}}\label{fig:vonkarman_chaotic}
\end{figure*}

We now extend the autonomous gSSM model to account for periodic forcing. A periodic force on the middle node of the beam, as shown in Fig. \ref{fig:vonkarman_chaotic}a, can make the dynamics of the beam chaotic. This has been observed and reported in the data-driven model of \cite{liuDatadrivenModelingForecasting2024}, who used simulation data of the full system. We now characterize the chaotic behavior without relying on any full-order simulations. 

To compare the full-order and the reduced-order dynamics, we take one of the buckled fixed points as an initial condition and simulate the beam under the influence of periodic forcing acting on the mid point with $\Omega = 25.3 \text{ rad}/s$ and $|\varepsilon \mathbf{f}_{ext}|=0.5 $N. Since this initial condition is known to be on the SSM, we run the same trajectory with the reduced dynamics after adding the leading-order contribution of the forcing as in \eqref{eq:forced_reducedDYn}.

The full-order trajectory is shown in Fig. \ref{fig:vonkarman_chaotic}b. As expected from the previous analysis and Fig. \ref{fig:vonkarman_buckled}, the forced SSM-reduced model blows up within a fraction of a second. In contrast, the globally valid gSSM-model exhibits sustained chaotic behavior, closely matching the full-order trajectory for short times. Fig. \ref{fig:vonkarman_chaotic}d-e shows that the chaotic attractor appearing as a result of the periodic forcing extends way outside the domain of convergence of the Taylor series of the SSM. 

We also construct the Poincaré-map of the full model and the gSSM-model by sampling the trajectories at multiples of the driving period $T = \frac{2\pi}{\Omega}$. Because the gSSM-model is a simple 2D ODE, we can sample the Poincaré-map with a fine resolution to obtain the structure in the reduced phase space shown in Fig. \ref{fig:vonkarman_chaotic}c. Overlaying the Poincaré-map obtained from the full system, we see a close correspondence with the predicted attractor. 

In addition, we estimate the leading Lyapunov exponents based on the exponential rate of divergence of initially close trajectories \cite{ott_2002} as $\lambda_{gSSM} = (3.0 \pm 0.02) \ 1/s$ and $\lambda_{full} = (3.1 \pm 0.05) \ 1/s$.  {The reduced model is in close agreement with the full model, even though the forced dynamics was approximated by simply projecting the forcing term onto the tangent space of the SSM according to \eqref{eq:forced_reducedDYn}. Incorporating higher-order corrections of the forced dynamics, as in \eqref{eq:phase_dep}, further improves the model, leading to more accurate short-time predictions. We present these comparisons, along with additional properties of the chaotic gSSM model, in the Supplementary Material.}

\subsection{Example 4: Data-driven model of an inverted flag experiment}\label{sec:invertedflag}

\begin{figure*}
    \centering
    \includegraphics[width = 0.99\linewidth]{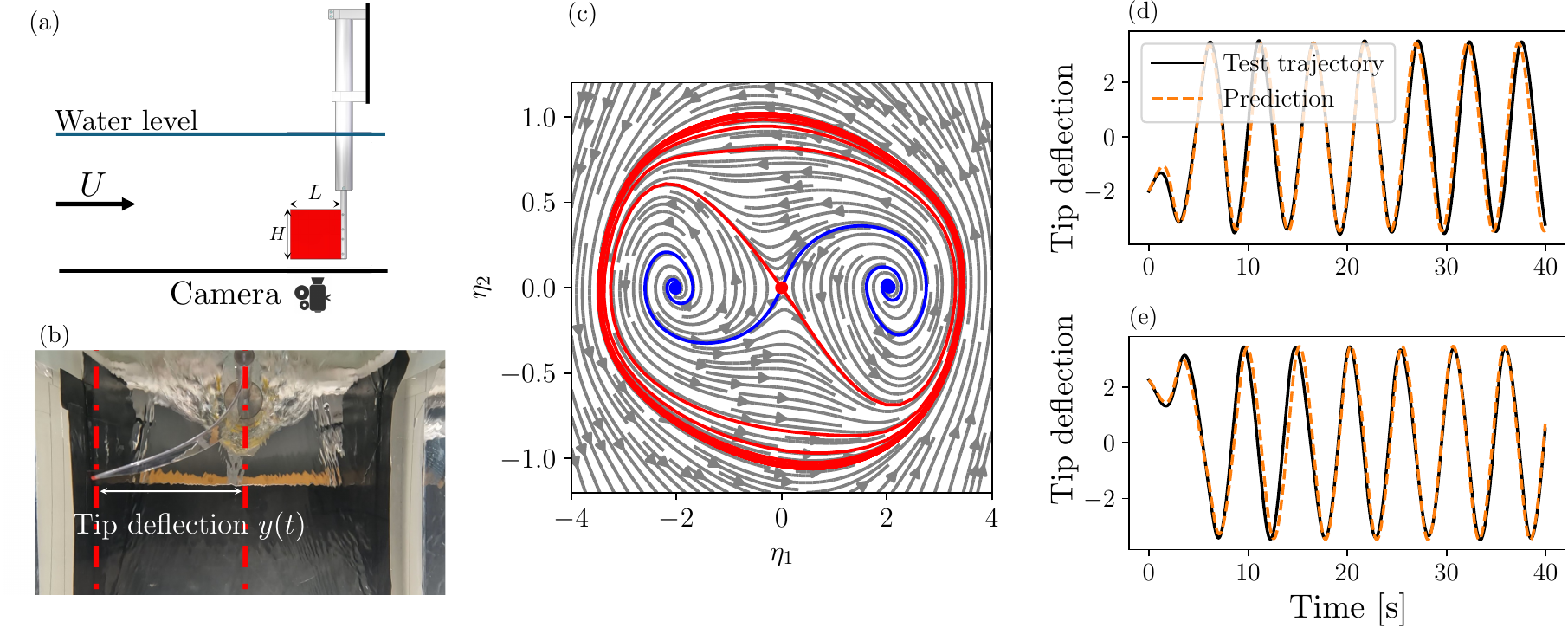}
    \caption{(a): Experimental setup of the inverted flag and snapshot of the experiment (b), courtesy of Giovanni Berti. The geometric parameters are $H=150$ mm, $L=150$ mm, $U=1\ \text{m}/\text{s}$. (c): Phase portrait of the gSSM-reduced dynamics of the inverted flag obtained from the $[5/5](\boldsymbol{\eta})$ approximation. The unstable fixed points are marked with colored dots. Blue curves denote the stable manifold of the undeflected state which connects to the two coexisting deflected fixed points. The red curve denotes the unstable manifold of the saddle, which wraps onto the stable limit cycle. (d)-(e): Predictions of the reduced model on test trajectories. The black curve is the true trajectory and the dotted orange curve is the gSSM-prediction.}\label{fig:invertedflag_reduced}
\end{figure*}

We now consider, as a data-driven example, the dynamics of an inverted flag, which is a flexible elastic sheet in the counterflow of a water tunnel. This configuration has generated recent interest due to its applications in energy harvesting \cite{windturbine} and vegetation \cite{vegetation}. From a modeling perspective, the inverted flag is a complicated fluid-structure interaction problem that benefits from reduced-order modeling. Here, we rely on experimental data obtained by \cite{xuDatadrivenModellingRegular2024}, whose experimental configuration is shown in Fig. \ref{fig:invertedflag_reduced}a-b. The elastic sheet is mounted in a water tunnel and is recorded from below with a video camera. Using classical image processing tools, the deflection of the tip of the flag $y(t)$ is recorded during the experiment. 

The main parameters of the system are the bending stiffness $K_B$ and the Reynolds number of the incoming flow, governed by the mean velocity $U$. By tracking the displacement of the tip of the flag, a data-driven SSM-reduced model was obtained \cite{xuDatadrivenModellingRegular2024} for both the large-amplitude periodic flapping regime and the chaotic flapping regime. In these regimes, the undeflected state of the flag is an unstable fixed point, which has a low-dimensional attracting mixed-mode SSM. 

We focus here on the large amplitude flapping regime with $K_B= 0.21$ and $\text{Re} = 6\times 10^4$. The slow SSM is tangent to a spectral subspace associated to an unstable real eigenvalue and a stable real eigenvalue and is two-dimensional. In addition to the saddle-type fixed point corresponding to the undeflected flag, the slow SSM contains two additional fixed points, which correspond to the deflected, but still stationary flag. The periodic flapping is a stable limit cycle within this slow SSM.

We now construct a data-driven gSSM-reduced model based on the tip-deflection data of \cite{xuDatadrivenModellingRegular2024}. A total of 16 trajectories are used for training. To reconstruct the slow SSM, we embed the tip displacement data $y(t)$ trajectories using $p=25$ time delays to form the observable $\mathbf{y}\in \mathbb{R}^{25}$. As in \cite{xuDatadrivenModellingRegular2024}, we approximate the tangent space of the SSM at the fixed point using the  {two leading principal components} of the delay-embedded trajectories. Since we use a moderate number of time delays and a short delay time, a leading-order (linear) approximation for the SSM suffices, i.e., we let 

\begin{equation}
    \begin{pmatrix}\eta_1 \\ \eta_2 \end{pmatrix} = \boldsymbol{\eta} = \mathbf{V}^T \mathbf{y}, \quad \mathbf{y} = \mathbf{V} \boldsymbol{\eta}. 
\end{equation}

The reduced dynamics are now approximated using rational functions, as detailed in Section \ref{sec:regression}, by solving the minimization problem \eqref{eq:ratfunc_error}. We find that a $[5/5]$ approximant gives the optimal reconstruction error, as computed on a validation trajectory, which was not used in the training. The resulting reduced vector field and the predictions on test trajectories are shown in Fig. \ref{fig:invertedflag_reduced}c-e.  {A diagonal approximant also ensures that the reduced dynamics remains well-behaved outside the domain of the training data. As opposed to classical polynomial regression, the values of the approximants remain bounded or only grow mildly for large $|\boldsymbol{\eta}|$. }

To obtain a comparable test accuracy to that shown in Fig. \ref{fig:invertedflag_reduced}d-e, previously, an order-11 polynomial approximation was used for the reduced dynamics. This required determining a total of $154$ coefficients $\mathbf{R}_{\mathbf{k}}$.  {The data-driven rational approximation, on the other hand, requires only $60$ coefficients, which is a significant reduction, for the cost of a slightly increased computational burden. In addition, the rational approximants extrapolate to larger domains in a more controlled way. } 

 {As we show in the Supplementary Material, preventing the rational functions from becoming singular is essential for an accurate approximation.} Other available methods, such as the rational function extension of the Sparse Identification of Nonlinear Dynamics (SINDy) algorithm \cite{manganInferringBiologicalNetworks2016a,bruntonDiscoveringGoverningEquations2016a} enforce no such constraints.

\section{Discussion}
We have presented a method to extend the range of validity of invariant manifold-based reduced-order models by applying Padé approximation, a classic analytic continuation technique. The Taylor coefficients of the parametrization of the SSM and the reduced dynamics were obtained with the robust numerical routines of \texttt{SSMTool}. We then extend their range of validity globally to obtain the gSSM-reduced model, by applying Padé approximation to these mappings that agree with the local Taylor expansions up to a prescribed order. 

We have demonstrated the method on high-dimensional examples of dynamical systems exhibiting global nonlinear behavior, such as transitions between steady states, large-amplitude oscillations, or even chaotic behavior. In all cases, gSSM-models obtained with Padé approximation significantly extend the domain of validity of the reduced model, reaching well beyond the domain of convergence of the classical Taylor series. We have also shown that the sign pattern of the univariate Taylor coefficients can be used to infer the location of the singularity.

To apply our approach, unexpected singularities in the approximants must be checked, and the optimal approximant of the Padé table must be determined according to the steps we have outlined in Section \ref{sec:pade}. Luckily, the numerical effort required to compute Padé approximants is negligible compared to the computation of the original Taylor coefficients of the SSM and the reduced dynamics. Therefore, even if no singularity-free Padé approximants can be found up to a given order, the benefits of a global SSM model far outweigh the costs of deriving it. 

We have shown that a data-driven analog of Padé approximation, rational function regression, can be directly applied to experimental data. Although rational function regression requires more computational power than the standard polynomial approximation implemented in \texttt{SSMLearn}, it can produce more accurate reduced models with fewer unknown coefficients, resulting in ultimately simpler models. 

We have demonstrated the method on the inverted flag example of \cite{xuDatadrivenModellingRegular2024} where the rational approximants of the SSM and the reduced dynamics were able to predict the coexisting deflected and undeflected fixed points and the stable limit cycle. As a main benefit, we note that the data-driven approach does not suffer from the same singularity issues as the equation-driven one. 

As we have demonstrated, Padé approximants and their data-driven extensions are particularly well-suited for approximating the global reduced dynamics on invariant manifolds. In addition to being able to represent more complex functions, they retain physical interpretability. In contrast to modern machine learning methods, the coefficients of the gSSM-reduced dynamics can be directly related to the underlying physics of the system. Specifically, the coefficients of the rational function can be interpreted as nonlinear frequencies or damping rates, as we have shown in Section \ref{sec:vonkarman}.

Finally, we note that although we have focused here on SSM-reduced models, the globalization method applies to other perturbative methods for dynamical systems, such as Poincaré-Lindstedt series \cite{falcoliniNumericalCalculationDomains1992}, geometric singular perturbation theory \cite{fenichelPersistenceSmoothnessInvariant1971,kuehnMultipleTimeScale2015}, or model reduction based on local linearization results \cite{haller_kaszas2024}.

\section{Methods}\label{sec:methods}
\subsection{Spectral submanifolds}\label{sec:ssms}
Let us assume that $\mathbf{x}=0$ is a hyperbolic fixed point of the system \eqref{eq:dynsys} with $\varepsilon=0$ and $\lambda_1, ..., \lambda_n\in \mathbb{C}$ are the eigenvalues of $\mathbf{A}$ with the corresponding eigenvectors denoted $\mathbf{e}_1, ..., \mathbf{e}_n\in\mathbb{C}^n$. We assume that a spectral gap condition holds for some $d<n$, i.e.,
\begin{equation}
    \text{Re }\lambda_n \leq ... \leq \text{Re }\lambda_{d+1} < \text{Re }\lambda_{d}\leq ... \leq \text{Re }\lambda_{1}.
\end{equation}
Let us denote the $d-$dimensional slow spectral subspace of $\mathbf{A}$ as $E$, which is defined as the span of the real and imaginary parts of the eigenvectors corresponding to the $d$ eigenvalues closest to zero. If the nonresonance condition 

\begin{align}
    &\lambda_k \neq \sum_{j=1}^n m_j \lambda_j, \quad k=1, ..., n, \\ 
   &m_j \in \mathbb{N}, \quad \sum_{j=1}^n m_j > 1, \nonumber
\end{align}
holds, a family of spectral submanifolds exist tangent to $E$, as discussed in  \cite{hallerNonlinearNormalModes2016b,hallerNonlinearModelReduction2023a}. The primary spectral submanifold is the unique, smoothest member of the family of $d-$dimensional invariant manifolds tangent to $E$ at the origin and is denoted as $\mathcal{W}(E)$ \cite{hallerNonlinearNormalModes2016b,hallerNonlinearModelReduction2023a}. 

Spectral submanifolds also exist for the non-autonomous system, i.e., for $\varepsilon>0$. In that case, SSMs are slow invariant manifolds attached to an anchor trajectory \cite{hallerNonlinearModelReduction2024}. For simplicity, we focus here on the case of periodic forcing with a single harmonic, where the anchor trajectory is a periodic orbit, and the original formulation of \cite{jainHowComputeInvariant2022b} applies. Denoting the parametrization of the slow SSM as $\mathbf{W}_\varepsilon(\mathbf{p}, \Phi)$, where $\Phi$ is the phase of the periodic forcing, the invariance equation reads as
\begin{align} \centering
    \label{eq:inveq_periodic}
     &\mathbf{A}\mathbf{W}_\varepsilon(\mathbf{p}, \Phi) + \mathbf{f}(\mathbf{W}_\varepsilon(\mathbf{p}, \Phi)) + \varepsilon \mathbf{f}_{\text{ext}}\cos \Phi  \nonumber \\ 
    &=D\mathbf{W}_\varepsilon(\mathbf{p}, \Phi)\dot{\mathbf{p}} + D_\Phi \mathbf{W}_\varepsilon \dot{\Phi},
\end{align}
 where $\dot{\Phi}=\Omega$ is the forcing frequency. The reduced dynamics is denoted as $\dot{\mathbf{p}} = \mathbf{R}_\varepsilon(\mathbf{p}, \Phi)$.

Due to the guaranteed smoothness of primary SSMs, the invariance equation can be solved using a Taylor expansion in $\mathbf{p}$ and a Fourier expansion in $\Phi$. The coefficients are then obtained by imposing the invariance equation order-by-order in the reduced coordinates.  {As shown by \cite{breunungExplicitBackboneCurves2018, ponsioenExactNonlinearModel2019, thurnherNonautonomousSpectralSubmanifolds2024}, to first order in $\varepsilon$ and to order-$\hat{N}$ in the reduced coordinates, the resulting expression is }
\begin{equation}
\label{eq:phase_dep}
    \mathbf{R}_\varepsilon(\mathbf{p}, \Phi) = \mathbf{R}(\mathbf{p}) +\varepsilon \sum_{|\mathbf{k}|=0}^{\hat{N}}\mathbf{S}_{\mathbf{k}}(\Phi) \mathbf{p}^\mathbf{k} + O(\varepsilon^2),
\end{equation}
 {where the coefficients $\mathbf{S}_{\mathbf{k}}(\Phi)$ can also be Fourier-expanded. Accounting for the phase-dependence of the nonautonomous coefficients in the reduced dynamics increases the accuracy of the reduced models for higher forcing amplitudes, at the cost of computing the coefficients $\mathbf{S}_{\mathbf{k}}(\Phi)$ for each forcing frequency of interest. These computations are already implemented in \texttt{SSMTool} \cite{thurnherNonautonomousSpectralSubmanifolds2024}.}

 {However, the phase-dependence of the parametrization and the reduced dynamics are often only small effects. It is, therefore, common to keep only the leading-order phase-dependent term in the reduced dynamics, resulting in 
\begin{align}
\label{eq:forced_reducedDYn}
   \mathbf{R}_\varepsilon(\mathbf{p}) = \mathbf{R}(\mathbf{p}) +\varepsilon\mathbf{V^*}\mathbf{f}_{ext}\cos (\Omega t),
\end{align}
where the operator $\mathbf{V}^*$ projects onto the tangent space and can be computed using the left eigenvectors of $\mathbf{A}$. Specifically, for a non-resonant two-dimensional SSM, we have 
\begin{align}
    \dot{\rho} &= \kappa(\rho)\rho + \varepsilon f \sin \psi \\
    \dot{\psi} &= \omega(\rho) - \Omega + \varepsilon f \frac{1}{\rho} \cos \psi,
\end{align}
where we have introduced the phase lag $\psi = \theta - \Omega t$ and the forcing amplitude $f$ is the projection of the external forcing amplitude \cite{jainHowComputeInvariant2022b}. }

The forced response of the system is then obtained by seeking fixed points of the reduced dynamics, which are given by the solutions of the equation $\dot{\rho} = \dot{\psi} = 0$. The amplitude $\rho_*$ of the response satisfies the implicit equation \cite{jainHowComputeInvariant2022b}
\begin{equation}
    \label{eq:forced_resp_eq}
    \left(\Omega - \omega(\rho_*)\right)^2 - \frac{\varepsilon^2 f^2}{\rho_*^2} + \kappa(\rho_*)^2 = 0.
\end{equation}

\subsection{Padé approximants}\label{sec:pade}
We now seek to improve the convergence properties of Taylor series approximations by summing the series outside the domain of convergence. One of the most popular methods of analytic continuation is Padé approximation, which is often used to sum divergent perturbative series. Padé approximation is primarily carried out on functions of a single complex variable, and hence it is directly applicable to 1D SSMs and to the functions $\kappa(\rho)$ and $\omega(\rho)$. 

Consider the function 
\begin{equation}
f(z): \mathbb{C} \to \mathbb{C},
\end{equation}
which could represent the parametrization of an invariant manifold, the reduced dynamics or $\kappa(\rho)$ and $\omega(\rho)$ in \eqref{eq:forced_resp_eq}. Let us denote its Taylor series representation around $z=0$ as
\begin{equation}
    f(z) = \sum_{n=0}^\infty c_n z^n.
\end{equation}
The Padé approximant of  {type $(N,M)$} is defined as the rational function 
\begin{equation}
    \label{eq:1d_pade_def}
    [N/M](z) = \frac{\sum_{n=0}^N a_n z^n}{\sum_{n=0}^M b_n z^n},
    \end{equation}
where $a_0=f(0)$ and $b_0$ can be chosen as $b_0=1$ without loss of generality. The Padé approximant is the best rational approximation of $f$ around $0$ \cite{bakerPadeApproximants1996} in the sense that its Taylor expansion around $0$ is the same as that of $f$ up to order $N+M$, i.e.,
\begin{equation}
    f(z)\sum_{n=0}^M b_n z^n - \sum_{n=0}^N a_n z^n = O(z^{N+M+1}).
\end{equation}
Therefore, the coefficients $b_n$ and $a_n$ can be determined by requiring
\begin{equation}
    \sum_{n=0}^\infty c_n z^n \sum_{m=0}^M b_m z^m = \sum_{n=0}^N a_n z^n
\end{equation}
for orders up to $N+M$ in $z$. This results in the linear equations for the coefficients $a_n$, $b_n$,
\begin{equation}
    \label{eq:1d_padecoeffs}
   a_n = \sum_{m=0}^n c_{n-m} b_{m}, \quad n=0,1,...,N+M.
    \end{equation}

We follow the robust approach of \cite{gonnetRobustPadeApproximation2013}, who first solve the homogeneous equations
\begin{equation}
    \label{eq:univariate_denom}
    \sum_{m=0}^n c_{n-m} b_{m} = 0, \quad n=N+1,...,N+M,
\end{equation}
 {with the convention of $b_j=0$ for $j<0$ or $j>M$. Using the singular value decomposition (SVD) of the Toeplitz matrix with elements [$c_{n-m}$], the coefficients $b_n$ can be computed. The remaining unknown coefficients $a_n$ are then given as
\begin{equation}
    \label{eq:univariate_num}
    a_n = \sum_{m=0}^n c_{n-m} b_{m}, \quad n=0,1,...,N.
\end{equation}}
This method is robust against numerical errors in the Taylor coefficients. In addition, many of the spurious poles of the approximant are removed, although not all of them \cite{mascarenhasRobustPadeApproximants2015}.

To apply Padé approximation to higher-dimensional SSMs, we need to generalize the method to multivariate power series. The multivariate generalization is not as straightforward as the univariate case, but multiple definitions exist. The most commonly used are {Chisholm approximants} \cite{chisholmRationalApproximantsDefined1974} and the homogeneous approximants \cite{bakerPadeApproximants1996, goncarLocalConditionSinglavaluedness1974, cuytHowWellCan1999}. Let us consider the multivariate function 

 \begin{equation}
    f: \mathbb{C}^d \to \mathbb{C}.
\end{equation}
given as a convergent Taylor series
\begin{equation}
    \label{eq:multi_f}
    f(\mathbf{z}) = \sum_{|\mathbf{k}|=0}^\infty c_{\mathbf{k}} \mathbf{z}^{\mathbf{k}}.
\end{equation}
We adopt the homogeneous approximants \eqref{eq:2dpade}, as defined originally by \cite{goncarLocalConditionSinglavaluedness1974, cuytHowWellCan1999} and require that the Taylor expansion of $[N/M](\mathbf{z})$ around $\mathbf{z}=0$ coincides with that of $f(\mathbf{z})$ up to order $N+M$, as 
\begin{align}
    \label{eq:multivariate_eq}
    f(\mathbf{z})\sum_{|\mathbf{k}|=0}^M b_{\mathbf{k}}\mathbf{z}^\mathbf{k} - \sum_{|\mathbf{k}|=0}^N a_{\mathbf{k}}\mathbf{z}^\mathbf{k} = O(\mathbf{z}^{N+M+1}) \\ \label{eq:multivariate_eq2}
    \sum_{|\boldsymbol{\ell}|=0}^{|\mathbf{k}|} c_{\mathbf{k}-\boldsymbol{\ell}}b_{\boldsymbol{\ell}} = a_{\mathbf{k}} \quad \text{ for } |\mathbf{k}|=0,1,...,N+M.
\end{align}

Since two-dimensional, oscillatory SSMs are the relevant objects for most systems, we consider the bivariate case as an illustration, i.e., with $\mathbf{z}=(z_1, z_2)^T$, where the coefficients are indexed on a lattice $\mathcal{L}_{N_1, N_2}$ defined as $\mathcal{L}_{N_1, N_2} = \{(i,j) \in \mathbb{N}^2: N_1 \leq i+j \leq N_2\} \subset \mathbb{N}^2$.

The conditions \eqref{eq:multivariate_eq2} become 
\begin{align}
    \sum_{k=0}^\alpha \sum_{\ell=0}^kc_{\beta-\ell, \alpha-k-\beta + \ell} b_{\ell, k-\ell} = a_{\beta, \alpha - \beta},
    \label{eq:2dpade_sols_num}
\end{align}
for $(\alpha,\beta) \in \mathcal{L}_{0, N}$ and 
\begin{align}
    \sum_{k=0}^\alpha \sum_{\ell=0}^kc_{\beta-\ell, \alpha-k-\beta + \ell} b_{\ell, k-\ell} = 0,
    \label{eq:2dpade_sols_denom}
\end{align}
for $(\alpha,\beta) \in \mathcal{L}_{N, N+M}$. Note, however, that the total number of unknowns and equations in \eqref{eq:2dpade_sols_num} and \eqref{eq:2dpade_sols_denom} are not equal. Therefore, instead, we seek a least-squares solution to \eqref{eq:2dpade_sols_denom} \cite{guillaumeMultivariatePadeAApproximation2000a}.

To summarize, gSSM-reduction consists of the following steps to compute the $[N/M]$ Padé approximants of the parametrization and the reduced dynamics of a slow SSM. 
\begin{enumerate}
    \item Compute the Taylor-series expansion of the SSM up to order $N+M$ using either the normal form style parametrization or the graph style parametrization. This returns $\mathbf{W}(\mathbf{p})$, $\omega(\rho)$ and $\kappa(\rho)$. 
    \item Solve the homogeneous linear equations \eqref{eq:2dpade_sols_denom} and \eqref{eq:univariate_denom} for the coefficients of the denominators of the parametrization and of $\omega$ and $\kappa$ using the robust SVD-based method \cite{gonnetRobustPadeApproximation2013}. 
    \item Check the zero sets of the denominator functions. If they contain points in the region of interest near the origin, adjust the orders $N$ and $M$. We found that diagonal and near diagonal $[N \pm 2/N\pm 2]$ approximants tend to work the best.
    \item Compute the coefficients of the numerators of the parametrization and of $\omega$ and $\kappa$ by evaluating \eqref{eq:2dpade_sols_num} and \eqref{eq:univariate_num}.  
\end{enumerate}

Although Padé approximation is a well-researched topic, convergence results are generally limited. For meromorphic functions $f(z)$ with a finite number of poles at $z_1, ..., z_k$, the theorem of de Montessus de Ballore \cite{demontessusdeballoreFractionsContinuesAlgebriques1905} guarantees the convergence of the approximants 
\begin{equation}
    [M/k](z) \text{ as } M \to \infty
\end{equation}
globally. A stronger result is available for Stieltjes functions, i.e. for functions of the form 
\begin{equation}
    f(z) = \int_{0}^\infty \frac{d\mu(t)}{1+zt},
\end{equation} for some positive measure $\mu$. In this case, the Padé approximants converge to $f(z)$ for all $z$ outside the negative real axis \cite{bakerPadeApproximants1996}. 

In the Supplementary Material, we discuss the Stiltjes-type center manifold of Euler's system  \cite{eulerSeriebusDivergentibus1760, bromwichIntroductionTheoryInfinite1947, arnoldGeometricalMethodsTheory1988, hallerNonlinearNormalModes2016b}, which, although non-analytic, can be described globally using Padé approximants.

\subsection{Rational function regression}\label{sec:regression}
 Let us denote the function to be approximated as $\mathbf{f}:\mathbb{R}^d \to \mathbb{R}^l$. This could represent $\mathbf{f}=\mathbf{W}$ with $l=n$ or $\mathbf{f} = \mathbf{R}$ with $l=d$.  {We assume that the value of $\mathbf{f}$, denoted $\boldsymbol{\zeta}_i = \mathbf{f}(\boldsymbol{\eta}_i)$, is known at points $\boldsymbol{\eta}_i, ..., \boldsymbol{\eta}_K$ in the domain of interest. We then approximate $\mathbf{f}$ as}
 \begin{equation}
    \label{eq:ratfunc_approx}
    \mathbf{f}(\boldsymbol{\eta}) \approx [N/M](\boldsymbol{\eta}) = \frac{\sum_{|\mathbf{k}|=0}^N \mathbf{a}_{\mathbf{k}}\boldsymbol{\eta}^{\mathbf{k}}}{\sum_{|\mathbf{k}|=0}^M {b}_{\mathbf{k}}\boldsymbol{\eta}^{\mathbf{k}}},
\end{equation}
where we have chosen a common denominator with coefficients $b_{\mathbf{k}}$ for all components of $\mathbf{f}$, similarly to vector Padé approximants \cite{sidiVectorExtrapolationMethods2017, laneSteadyBoussinesqConvection2024}. To avoid introducing singularities for the approximants, we require that the denominator is non-zero at all points $\boldsymbol{\eta}_i$. This is equivalent to requiring that the denominator is strictly positive. 

The coefficients are determined by minimizing the error 
\begin{align}
    \label{eq:ratfunc_error}
    \mathcal{E}_r &= \sum_{i=1}^K \left| \boldsymbol{\zeta}_i  - \frac{\sum_{|\mathbf{k}|=0}^N \mathbf{a}_{\mathbf{k}}\boldsymbol{\eta}_i^{\mathbf{k}}}{\sum_{|\mathbf{k}|=0}^M {b}_{\mathbf{k}}\boldsymbol{\eta}_i^{\mathbf{k}}} \right|^2, \text{ such that} \\
     & \sum_{|\mathbf{k}|=0}^M {b}_{\mathbf{k}}\boldsymbol{\eta}_i^{\mathbf{k}} \geq \delta \quad \text{ for } i=1,...,K,
\end{align}
for some small $\delta>0$. The constrained minimization problem is solved by gradient-based optimization methods \cite{batesNonlinearRegressionAnalysis1988}. The initial guess for the coefficients is obtained by solving the linearized problem, minimizing 
\begin{equation}
    \sum_{i=1}^K \left| \left(\sum_{|\mathbf{k}|=0}^M {b}_{\mathbf{k}}\boldsymbol{\eta}_i^{\mathbf{k}}\right)\boldsymbol{\zeta}_i  - \sum_{|\mathbf{k}|=0}^N \mathbf{a}_{\mathbf{k}}\boldsymbol{\eta}_i^{\mathbf{k}} \right|^2,
\end{equation}
subject to the same positivity constraint. This is a linear least-squares problem, which can be solved efficiently using the method of \cite{austinPracticalAlgorithmsMultivariate2021}, based on robust Padé approximation \cite{gonnetRobustPadeApproximation2013}.

\section{Acknowledgements}
We are grateful to Mohammad Farazmand for sharing his code to simulate the Kolmogorov flow.
\section{Data availability}
All data processed in the manuscript will be made available upon publication.
\section{Code availability}
The codes containing the implementation of the methods will be made available upon publication.

\newpage

{\bf \Huge \centering Supplementary Material}
\vspace{40pt}
\section{Further examples}

\subsection{The Euler example}\label{sec:Euler}
\begin{figure}[h!]
    \centering
    \includegraphics[width = 0.99 \linewidth]{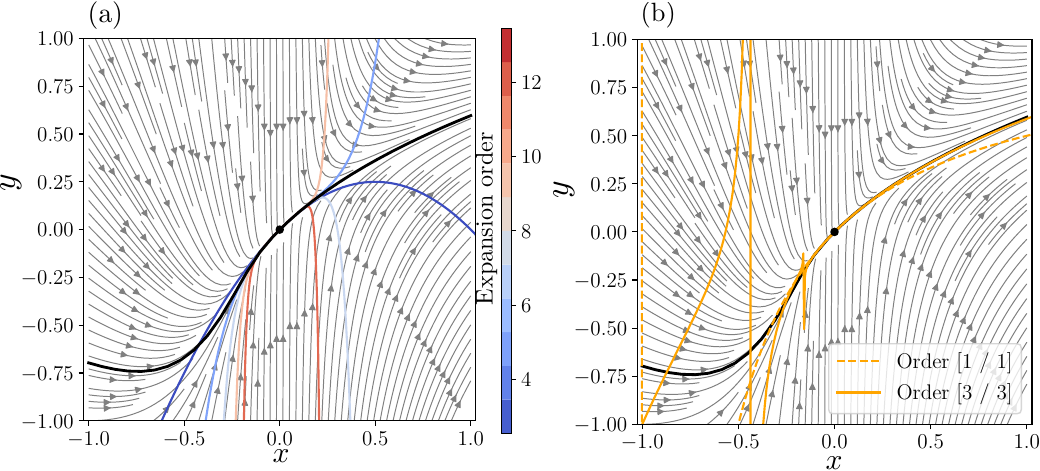}
    \caption{(a): Phase portrait of \eqref{eq:eulereq} with the $C=0$ member of the center manifold family \eqref{eq:euler_center} shown in black. (a): Evaluation of the Taylor expansion \eqref{eq:euler_mdf_exp} up to order 15. (b): Padé approximants $[1/1](x)$ and $[3/3](x)$ of the series \eqref{eq:euler_mdf_exp}.}
    \label{fig:euler_example}
\end{figure}

We now recall an example from \cite{arnoldGeometricalMethodsTheory1988} and \cite{hallerNonlinearNormalModes2016b}, which has originally been studied by Euler \cite{eulerSeriebusDivergentibus1760}. The system of equations reads as
\begin{align}\label{eq:eulereq}
    \dot{x} &= x^2 \\
    \dot{y} &= x - y.
\end{align}
The origin is a non-hyperbolic fixed point having a one-dimensional $C^\infty$ center manifold. In fact, the whole phase space is foliated by one-dimensional center manifolds, and there is no distinguished one in terms of smoothness \cite{hallerNonlinearNormalModes2016b, vanstrienCenterManifoldsAre1979}. Nevertheless, the center manifold can be parametrized as $y=h(x)$, which leads to the invariance equation
\begin{equation}
\label{eq:inveq_euler}
    h(x)'x^2 = x - h(x), 
\end{equation}
which was studied by \cite{eulerSeriebusDivergentibus1760}. 
Seeking a power series approximation of $h(x)$ leads to 
\begin{equation}
\label{eq:euler_mdf_exp}
    y=h(x)\sim\sum_{k=0}^\infty (-1)^k (k-1)! x^k \text{ as }x\to 0,
\end{equation}
where the notation $\sim$ means that the asymptotic series need not converge for $x\neq 0$. The invariance equation \eqref{eq:inveq_euler} can be solved, for example, by multiplying by the integrating factor $e^{-1/x}$, which leads to the expression
\begin{equation}
    \label{eq:euler_center}
    h(x,C) = Ce^{\frac{1}{x}} - e^{\frac{1}{x}}\text{Ei}\left(-\frac{1}{x} \right) 
\end{equation}
where $C\in \mathbb{R}$ is arbitrary and
\begin{equation}
    \text{Ei}(z) = \int_{-\infty}^x \frac{e^{t}}{t}dt
\end{equation}
denotes the exponential integral \cite{arfkenMathematicalMethodsPhysicists2005}. The center manifolds all have the coinciding asymptotic expansion \eqref{eq:euler_mdf_exp}, which has a zero radius of converge, as seen also in Fig. \ref{fig:euler_example}a. 

Euler derived the differential equation \eqref{eq:inveq_euler} to assign a value to the sum of factorials with alternating signs, i.e., to the sum $1 - 2! + 3! - 4! + ...$. We refer to Chapter XII in \cite{bromwichIntroductionTheoryInfinite1947} for historical context. 

Computing Padé approximants of the series \eqref{eq:euler_mdf_exp}, however, reveals that even low-order approximants approximate the $C=0$ member of the center manifold family for $x>0$. Indeed, $h(x,0)$ is an example of a Stiltjes function \cite{benderAdvancedMathematicalMethods1999} and hence belongs to the rare class of functions for which the convergence of Padé approximants can be proven. Padé approximants of Stiltjes functions converge everywhere in the complex plane except along the negative real axis, where the original function $h(x,0)$ has a branch cut. 

This behavior is shown in Fig. \ref{fig:euler_example}b, where we see that even the low-order approximants show good agreement for $x>0$. Their poles are concentrated on $x<0$, mimicking the branch cut of the original function $h(x,0)$. The approximants shown can be written as 

\begin{equation}
    y=[1/1](x) = \frac{x}{1+x} \quad y=[3/3](x) = \frac{x^3 + 8x^2+11x}{6x^3 +18x^2+9x+1}
\end{equation}

\subsection{The Dauchot-Manneville model}

The Dauchot-Manneville model is a planar model given by the system of equations
\begin{equation}
    \label{eq:dauch_manneville}
    \begin{pmatrix} \dot{x}_1 \\ \dot{x}_2 \end{pmatrix} = \begin{pmatrix}s_1 & 1 \\ 0 & s_2 \end{pmatrix}\begin{pmatrix} x_1 \\ x_2 \end{pmatrix} + \begin{pmatrix}
        x_1 x_2 \\ -x^2_1
    \end{pmatrix}.
\end{equation}
The model is bistable and shows strong non-normality at the trivial fixed point \cite{dauchotLocalGlobalConcepts1997}. These properties  are reminiscent of the subcritical transition to turbulence in shear flows \cite{halcrowHeteroclinicConnectionsPlane2009}. It is 
For the parameter values $s_2<0, s_1<0$, this system has two coexisting stable fixed points, $p_1=(0,0)$ and $p_3$, and a saddle point, $p_2$. A one-dimensional manifold tangent to the slow spectral subspace connects the saddle to the two stable fixed points. We now compute the slow SSM of the origin, which is the same slow manifold observed by \cite{dauchotLocalGlobalConcepts1997}, hoping that this coincides with the connecting orbit. The nonresonance conditions of \cite{hallerNonlinearNormalModes2016b} are satisfied, and hence the slow SSM is analytic near the origin.

To find the parametrization of the slow SSM we assume that it can be written as a graph over the variable $x_1$. To construct this parametrization we assume a Taylor expansion of the form

\begin{equation}
    \label{taylor_DM_model}
    x_2 = W(x_1) = \sum_{n=1}^\infty w_n x_1^n
\end{equation}
and find the coefficients $w_n$. Substituting \eqref{taylor_DM_model} gives, to cubic order, \begin{equation}
    \label{taylor_DM_model_leading}
    h(x_1) = -\frac{x_{1}^{2}}{2 s_{1} - s_{2}} - \frac{2 x_{1}^{3}}{\left(2 s_{1} - s_{2}\right)^{2} \left(3 s_{1} - s_{2}\right)}  + O(x_1^4).
\end{equation}

\begin{figure}
    \centering
    \includegraphics[width = 0.99 \linewidth]{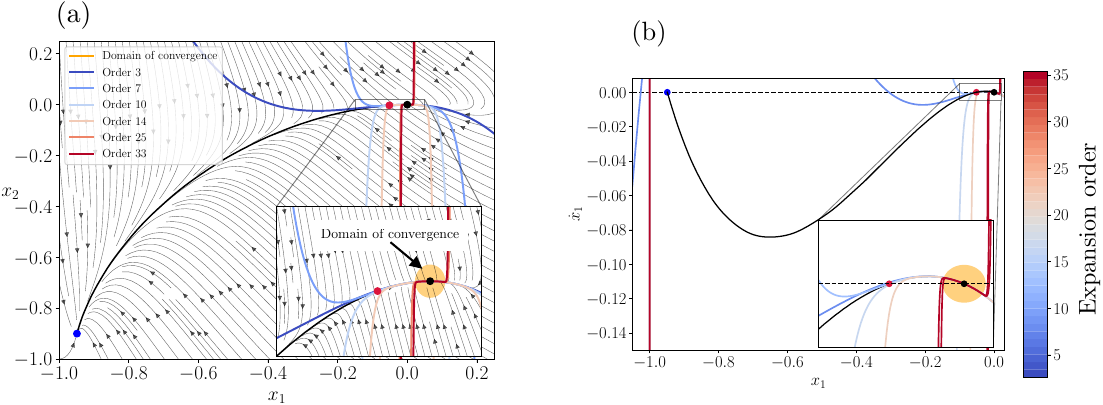}
    \caption{(a) Phase portrait of the Dauchot-Manneville model \eqref{eq:dauch_manneville} with $s_1=-0.038$, $s_2 = -1$ and the Taylor-series approximations of the slow SSM of the stable fixed point. The inset shows the domain on which the Taylor series converges. (b): The reduced dynamics on the SSM, obtained from the Taylor series approximation. The saddle is shown in red, the stable fixed points in blue, and black. The inset shows the same dynamics in the neighborhood of $p_1$ and $p_2$.}
    \label{fig:dauch_taylor}
\end{figure}

Fixing $s_1=-0.038$, in Fig. \ref{fig:dauch_taylor}, we show the manifold computed up to increasing orders in the Taylor expansion, which is guaranteed to have a nonzero radius of convergence. To estimate the radius of convergence, \cite{kaszasNewMethodsReducedorder2023} computed the accummulation point of zeros of the Taylor-expansion \eqref{taylor_DM_model} as in \cite{ponsioenModelReductionSpectral2020}. In addition, we find that all coefficients $w_n$ are negative. Therefore, by the Vivanti–Pringsheim theorem \cite{remmertTheoryComplexFunctions1998}, the convergence limiting singularity is at some small positive $x_1^s>0$, which is approximately given as $x_1^{s}=0.01$. This limits the convergence even for $x_1<-x_1^s$, and therefore, the reduced-order model cannot contain any of the non-trivial fixed points.

The reduced dynamics on the SSM is obtained by substituting $x_2 = h(x_1)$ into \eqref{eq:dauch_manneville}, which leads to 
\begin{equation} 
    \label{dauch_manneville_reduced}
\dot{x}_1 = s_1 x_1 + h(x_1) + x_1 h(x_1). 
\end{equation}
The roots of \eqref{dauch_manneville_reduced} correspond to the predicted fixed points along the SSM. We show the reduced dynamics in Fig. \ref{fig:dauch_taylor}b. As expected, outside the domain of convergence, the dynamics predicted from the Taylor expansion is very different from the true reduced dynamics computed numerically. 

We now to re-sum the divergent Taylor series outside its domain of convergence by computing its Padé approximants given as \eqref{eq:1d_pade_def}. Given the Taylor series representation of $h(x_1)$, the Padé approximants can be constructed by solving the linear equations \eqref{eq:1d_padecoeffs}. The sequence of the diagonal approximants, i.e. when $N=M$ generally has better convergence and approximation properties \cite{benderAdvancedMathematicalMethods1999}.
\begin{figure}[h!]
    \centering
    \includegraphics[width = 0.99 \linewidth]{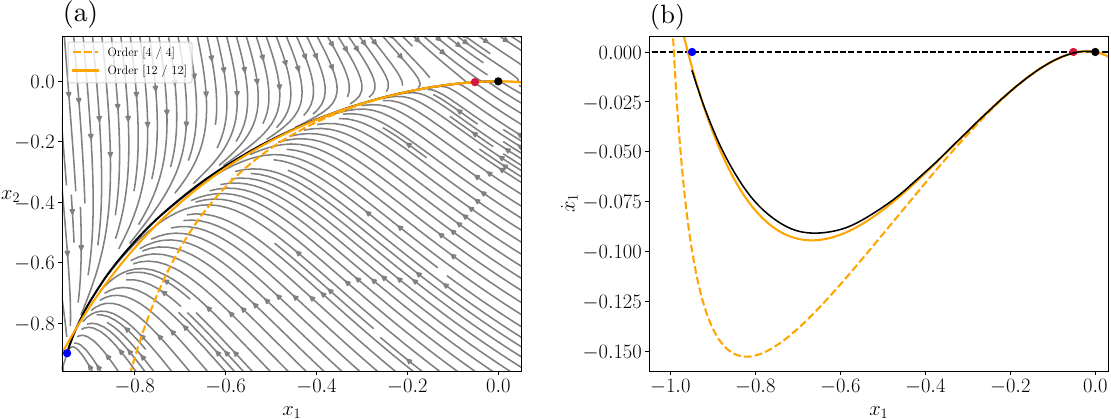}
    \caption{Padé approximations of the slow SSM in the Dauchot-Manneville model with the same parameters as in Fig. \ref{fig:dauch_taylor}. (a): the approximation in the phase space. (b): the reduced dynamics. }
    \label{fig:dauchot_pade}
\end{figure}
The diagonal Padé approximants of orders $M=N=4$ and $M=N=12$ are shown in Fig. \ref{fig:dauchot_pade}. The lower-order approximant can be written as 
\begin{equation}
    h(x_1) = \frac{-1.1x_1^2 + 40.55 x_1^3 - 295.25 x_1^4 }{1-39.2x_1 + 347.58 x_1^2 - 324.01 x_1^3 - 647.28 x_1^4}.
\end{equation}

The approximation of the manifold is remarkably good over the entire region of interest spanning all three fixed points for a high-order approximant. Note that even the low-order approximant correctly predicts the three fixed points as well as their stability types, albeit with a non-negligible error in their location. Thus we successfully extended the local information contained in the Taylor-expansion of the SSM to achieve a globally valid representation.

\subsection{Invariant manifold with an imaginary singularity}\label{sec:imsing}
It is often the case, that the convergence limiting singularity is not on the real axis, which means it is non-physical. As an example, consider the system 

\begin{align}
\label{eq:complex_sing}
    \dot{x} &= x \\ 
    \dot{y} &= -y + \frac{2x}{(x^2+1)^2}
\end{align}
It can be verified that the unstable manifold of the origin is parametrized as 
\begin{equation}
    y = h(x) = \frac{x}{1+x^2} = \sum_{n=0}^\infty (-1)^{2n}x^{2n+1}
\end{equation}

Note, however, that the first equality is valid for all $x\in\mathbb{R}$, the infinite sum representation is only valid for $|x|<1$. Denoting the complex extension of the parametrization as $\hat{y}=h(\hat{x})$, it is apparent that the convergence of the series is limited by the two poles of $h(\hat{x})$ at $\hat{x}=\pm i$. Therefore, the radius of convergence is 1. There is no singularity along the real axis, and the parametrization remains well-defined for any $x\in \mathbb{R}$. 

\begin{figure}[h!]
    \centering
    \includegraphics[width = 0.99 \linewidth]{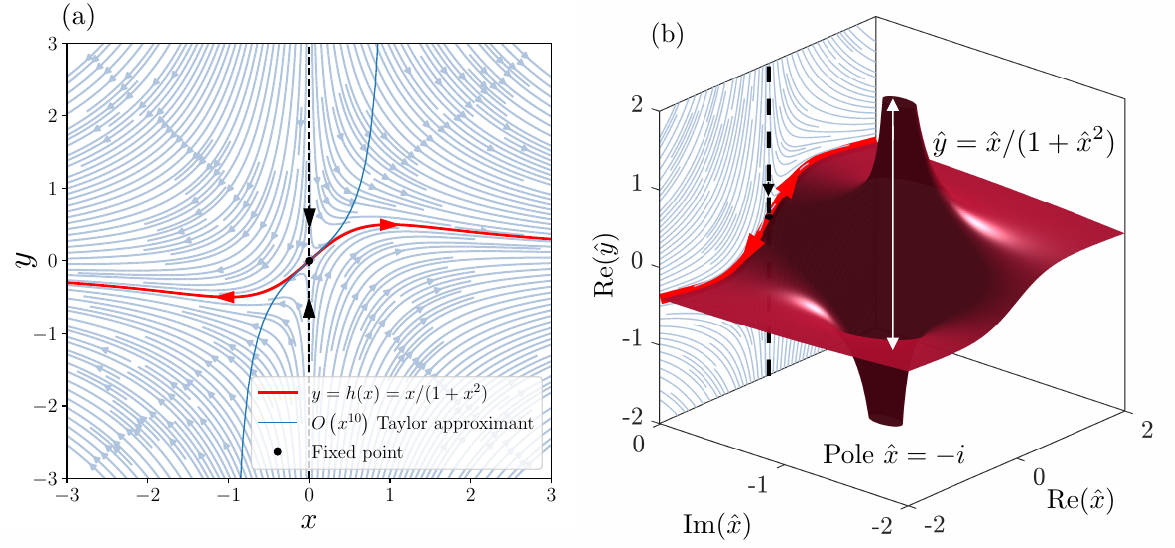}
    \caption{(a): Phase portrait of \eqref{eq:complex_sing} with the unstable manifold of the origin (red) and its order-10 Taylor approximant. (b): Visualization of the complex extension of the unstable manifold and the convergence limiting singularity at $\hat{x}=-i$. }
    \label{fig:complex_sing}
\end{figure}

As expected, the Taylor-series of $h(x)$ converges slowly, and only for $|x|<1$. This is shown in Fig. \ref{fig:complex_sing}. On the other hand, any Padé approximant $[N/M]$ with $N\geq1$ and $M\geq 2$ is exact and can be evaluated for any $x\in \mathbb{R}$. 

 \subsection{Shaw-Pierre oscillator} \label{sec:shawpierre}

\begin{figure}[H]
    \centering
    \includegraphics[width = 0.99 \textwidth]{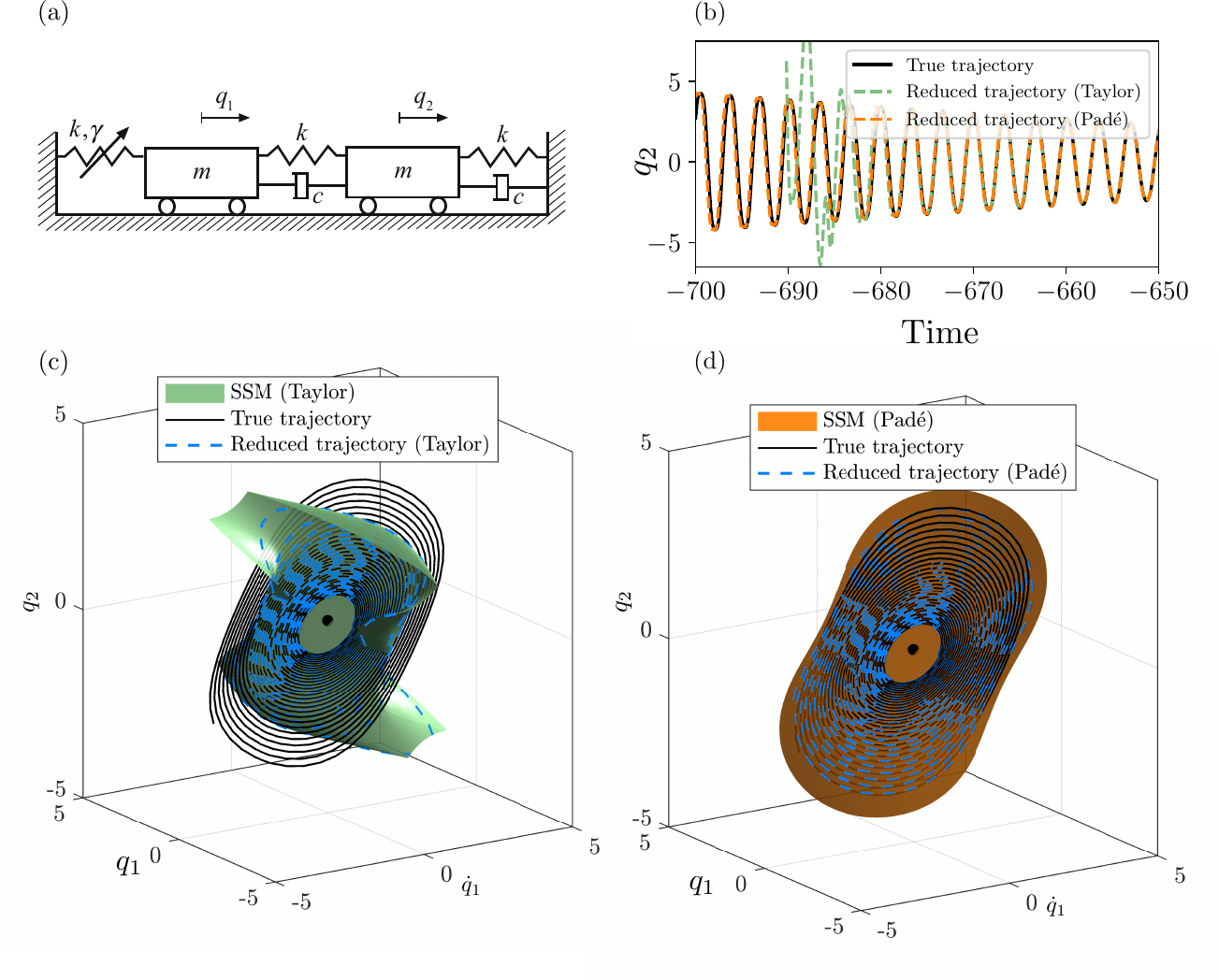}
    \caption{(a): Schematic diagram of the Shaw-Pierre system. (b): Time series of a backward integrated trajectory started from $p_0 = \bar{p_0} = 10^{-3}$ in the full model (black), in an order-18 Taylor approximation (green) and the corresponding [5/5] Padé approximant (orange). (c): The trajectory in the phase space and the Taylor-approximation of the SSM. (d): The Padé approximant of the SSM. }
    \label{fig:sharpierre_aut}
\end{figure}

We recall the Shaw-Pierre system from \cite{hallerNonlinearNormalModes2016b}, which is the forced and damped two-degree-of-freedom oscillator shown in Fig. \ref{fig:sharpierre_aut}a, whose equations of motion are given by
\begin{align}
    \label{eq:shawpierre}
    \ddot{q}_1 & + c(2\dot{q}_1 - \dot{q}_2) + k(2q_1 - q_2) + \gamma q_1^3 = \varepsilon \cos (\Omega t) \\
    \ddot{q}_2 & + c(2\dot{q}_2 - \dot{q}_1) + k(2q_2 - q_1)  =0,
\end{align}
with $k=3$, $c=0.003$, $\gamma = 0.5$. The phase space is spanned by the four variables $(q_1, \dot{q}_1, q_2, \dot{q}_2)$. For $\varepsilon =0$ the system has a stable fixed point at $(0,0,0,0)$. For $\varepsilon >0$ this fixed point perturbs into a periodic orbit, which represents the forced response of the system. Using SSMTool, we can calculate the two-dimensional SSM of the autonomous system ($\varepsilon=0$).

 To visualize the breakdown of the convergence of the Taylor series for high amplitudes $|p|$, we select an initial condition close to the origin of the unforced system at $p_0 = \bar{p_0} = 10^{-3}$ and integrate the system backwards in time under both the reduced- and the full-order dynamics. The resulting trajectory and its Taylor approximation are shown in Fig. \ref{fig:sharpierre_aut}b and Fig. \ref{fig:sharpierre_aut}c. The Taylor-approximation fails to describe an invariant manifold for $q_2>3$, which is also manifested by apparent self intersections of the manifold.  

\begin{figure}[H]
    \centering
    \includegraphics[width = 0.99\textwidth]{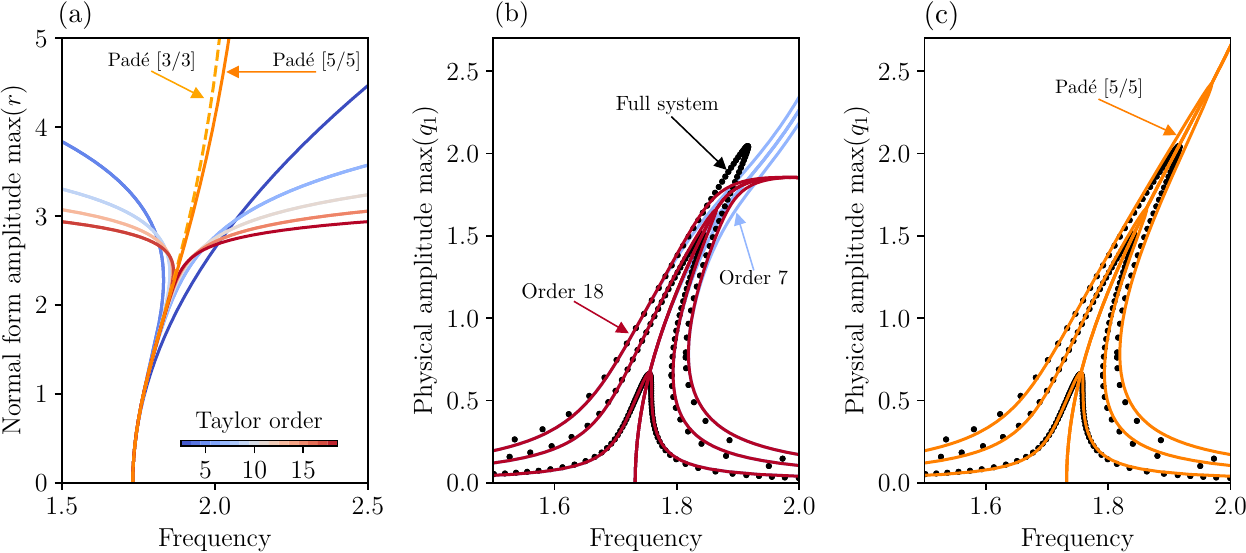}
    \caption{Padé approximants of the backbone curves and forced response curves of the Shaw-Pierre oscillator. (a): Backbone curves expressed as $\omega(\rho)$. (b): Full-order response (black) with $\varepsilon = 0.05, 0.2, 0.3$ and its order-7 and 18 Taylor approximation. (c): Padé approximants of the response. }\label{fig:shawpierre_backbone}
\end{figure}

The backbone curve, i. e. the function $\omega(\rho)$ is computed as 
\begin{equation}
    \label{eq:shawpierreomega}
    \omega(\rho) = 1.7320 + 0.0385 \rho^2 - 0.0037 \rho^4 + 0.0004 \rho^6 + O(\rho^8),
\end{equation}
and is shown in Fig. \ref{fig:shawpierre_backbone}a. The plots show that the Taylor series approximation of the backbone curve converges for amplitudes up to $\rho \approx 3$. Moreover, the sign pattern of \eqref{eq:shawpierreomega} shows alternating positive and negative coefficients, suggesting that the convergence limiting singularity is along the imaginary axis, similarly to the example in Section \ref{sec:imsing}. This can be inferred from the calculations of \cite{mercerCentreManifoldDescription1990b}, which we specialize to our case in Section \ref{sec:compute_sing}.

Computing the forced response curves for increasing orders of the Taylor-expansion, we find in Fig. \ref{fig:shawpierre_backbone}b that they agree with the full-order computations as long as the response amplitude is smaller than $q_1 \approx 2$, as also inferred from the autonomous analysis in Fig. \ref{fig:sharpierre_aut}.  The full order computations were carried out using the numerical continuation package COCO \cite{dankowiczRecipesContinuation2013}.

Padé approximants, however have no problem approximating the true forced response even outside this domain of convergence. We compare the [3/3] and [5/5] approximants, which closely agree with each other. The forced response curves are shown in Fig. \ref{fig:shawpierre_backbone}c. For example, the $[5/5]$ approximant of the function $\omega(\rho)$ is given as
\begin{equation}
    \label{eq:shawpierre_pade}
    [5/5](\rho) = \frac{1.7320 + 0.3717 \rho^2 + 0.0166 \rho^4}{1 + 0.1924 \rho^2 + 0.0074\rho^4}
\end{equation}

\section{Kolmogorov flow}
 The Kolmogorov flow, in Fourier space, is governed by the following ordinary differential equations
\begin{align}
    \label{eq:kolmogorov_fourier}
    \frac{d\hat{\omega}}{dt} = -\frac{1}{\text{Re}}(k_x^2 + k_y^2) \hat{\omega} - \widehat{(\mathbf{u}\cdot \nabla )\omega} - 4 \delta_{k_y, 4} \delta_{k_x,0},
\end{align}
for the discrete wave numbers $k_x,k_y=-12, ..., 12$, resulting in a total of $576$ degrees of freedom. Based on \cite{chandlerInvariantRecurrentSolutions2013, farazmandTensorbasedFlowReconstruction2023} we use a pseudo-spectral implementation. The nonlinearity is evaluated by introducing the stream function $\psi(x,y)$ as 
\begin{equation}
    \mathbf{u} = \begin{pmatrix} \partial_y \psi  \\ -\partial_x \psi \end{pmatrix}.
\end{equation}
The streamfunction can be recovered by solving the Poisson equation in Fourier space
\begin{equation}
     (k_x^2 + k_y^2)\hat{\psi} = -\hat{\omega}.
\end{equation}
The inverse Fourier transform of $\hat{\psi}$ then allows us to evaluate the nonlinearity $(\mathbf{u}\cdot \nabla )\omega$ in real space and take the Fourier transform subsequently. For dealiasing, we use the 3/2 scheme

\section{Locating the convergence limiting singularity}\label{sec:compute_sing}
The location of the convergence limiting singularity of a Taylor-series expansion can be inferred from the sign pattern and the size of the Taylor-coefficients. Here our main focus is to decide whether the singularity is on the real axis, as that would potentially indicate a genuine singularity of the function. This can be done by analyzing the sign pattern of the Taylor-coefficients.

We specialize the analysis of \cite{mercerCentreManifoldDescription1990b} to the case of Taylor-approximations of backbone curves and damping curves obtained from two-dimensional SSM-reduced models. The backbone curve $\omega(\rho)$ and the damping curve $\kappa(\rho)$ are even functions of their arguments. The sign pattern is governed by the closest singularity of the function in the complex plane. Consider the prototype even function with poles at $z=\pm re^{\pm i \theta}$, singularity given by 

\begin{align}
    G(z) &= g(z) + g(-z), \text{ with}  \nonumber\\ 
    g(z) &= \left(1 - \frac{z}{re^{i\theta}}\right)^\nu + \left(1- \frac{z}{re^{-i\theta}}\right)^\nu. 
\end{align}
The Taylor-expansion of $g(z)$ around $z=0$ converges for $|z|<r$ and is given as \cite{mercerCentreManifoldDescription1990b}
\begin{equation}
    g(z) = \sum_{n=0}^\infty 2(-1)^{n} \binom{\nu}{n}r^{-n}\cos(n\theta) z^n,
\end{equation}
which allows us to write
\begin{equation}
    G(z) = \sum_{n=0}^\infty 2 \binom{\nu}{2n}r^{-2n}\cos(2n\theta) z^{2n}.
\end{equation}

The sign pattern of the Taylor-coefficients is determined by the cosine term. In our examples, the signs of the Taylor-coefficients of the backbone curve and the damping curve of the von Kármán beam, and the Shaw-Pierre oscillator showed an alternating sequence of positive, zero, and negative values. This indicates that $\theta = \pi / 2$, i.e. the singularity is likely close to the imaginary axis.

\section{Checking for singularities}
\begin{figure}[h!]
    \centering
    \includegraphics[width = 0.99\textwidth]{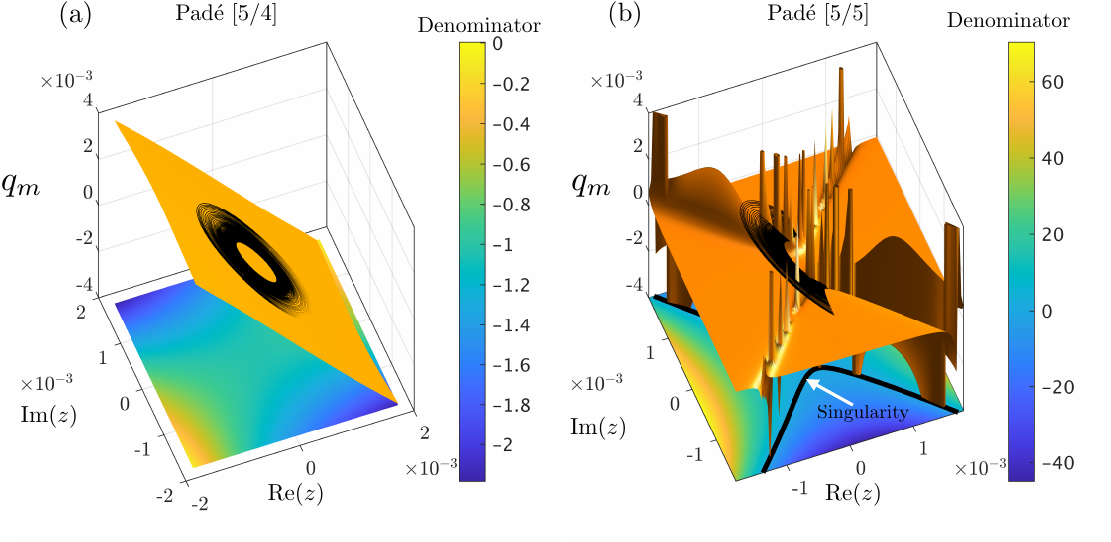}
    \caption{Invariance of the Padé approximant of the parametrization. The orange surface is the Padé approximant corresponding to the tip of the beam. The value of the denominator is also shown, color-coded. (a): the [5/4] approximant is well-behaved near the origin and the backward trajectory (black). (b): The [5/5] approximant has a spurious singularity corresponding to the zero set of the denominator. }\label{fig:vonkarman_singularity}
\end{figure}

Special care must be taken to avoid spurious singularities of Padé approximants, especially in the multivariate case. For the parametrization of the von Kármán beam, we chose the $[5/4]$ Padé approximant.
 This was because the diagonal approximant $[5/5]$ had a spurious curve of singularities in the neighborhood of the origin, as shown in Fig. \ref{fig:vonkarman_singularity}b. Since no general pointwise convergence result is available for the Padé approximants, one must always check whether singularities are present in the approximants. We found that decreasing the order of the denominator by one was sufficient to eliminate this singularity, as shown in Fig. \ref{fig:vonkarman_singularity}a. To verify the invariance of the manifold, we compute a backward-trajectory started very close to the origin, as shown in the main text. 

\section{Properties of the chaotic von Kármán beam}
To further analyze the chaotic properties of the reduced model on the gSSM arising due to forcing. Specifically,  we estimate the leading Lyapunov exponent of the chaotic attractor observed in the full model and the gSSM-reduced model. We perturb the initial conditions of the trajectories presented in the main text by a small perturbation of size $10^{-7}$ in the reduced coordinates. The rate at which the reference and the perturbed trajectories deviate is governed by the leading Lyapunov exponent of the attractor \cite{ott_2002}.

\begin{figure}[h!]
    \centering
    \includegraphics[width = 0.99\textwidth]{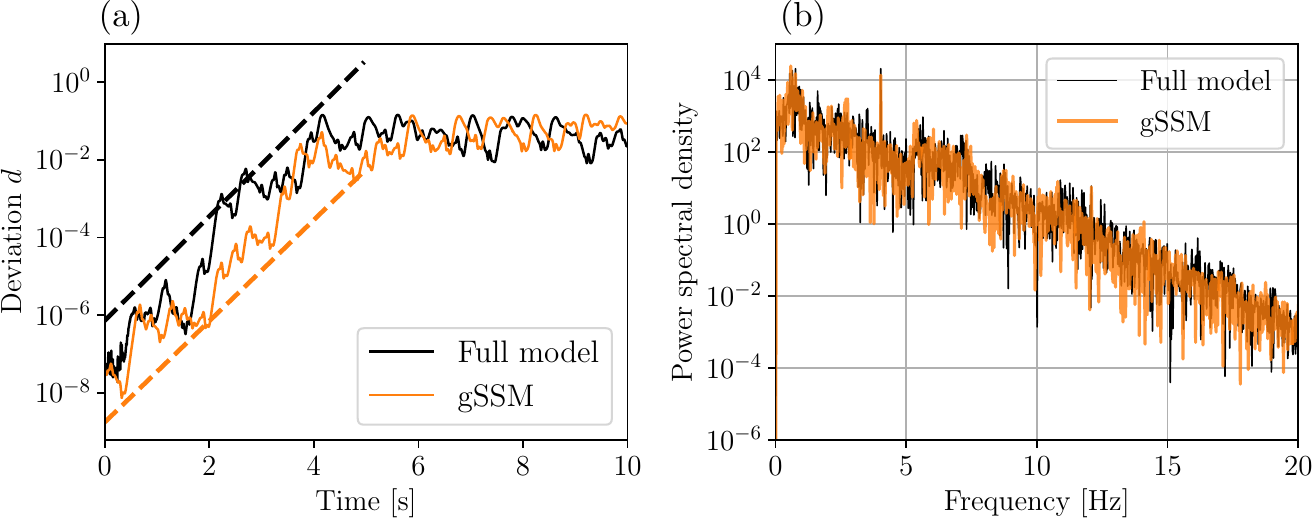}
    \caption{(a): Estimation of the leading Lyapunov exponent based on the rate of divergence of nearby trajectories. Exponents are fitted to the initial exponential trend of the curves. We obtain $\lambda_{gSSM} = (3.0 \pm 0.02) \ 1/s$ and $\lambda_{full} = (3.1 \pm 0.05) \ 1/s$. The corresponding exponential functions are indicated with dashed lines.  (b): Power spectral density of the attractors of the reduced and full-order models. }\label{fig:vonkarman_chaotic_props}
\end{figure}

Denoting the instantaneous distance between the perturbed and the reference trajectories as $d(t)$, we expect that $d\sim e^{\lambda t}$ holds with $\lambda>0$ on the attractors of both the full system and the gSSM-reduced system. 
Figure \ref{fig:vonkarman_chaotic_props}a shows the estimation of the leading Lyapunov exponent. 

In addition, we also estimate the power spectral density (PSD) of the chaotic attractors. Based on long simulations of the reduced and full-order models lasting 200 forcing periods, we compute the frequency spectrum of the $\eta_1$ component of this time series. The power contained in this spectrum is shown in Fig. \ref{fig:vonkarman_chaotic_props}b. The power spectrum computed on the gSSM closely matches the full-order simulation. The spectrum features a wide range of frequencies, indicative of chaotic behavior.

\subsection{Nonautonomous model}
The reduced dynamics on the SSM reads as 
\begin{equation}
\label{eq:phase_dep}
    \mathbf{R}_\varepsilon(\mathbf{p}, \Phi) = \mathbf{R}(\mathbf{p}) +\varepsilon \sum_{|\mathbf{k}|=0}^{\hat{N}}\mathbf{S}_{\mathbf{k}}(\Phi) \mathbf{p}^\mathbf{k} + O(\varepsilon^2),
\end{equation}
where $\hat{N}\geq 0$ is the approximation order for the forcing. In the main text, we have only considered the leading-order contribution with $\hat{N}=0$, but the accuracy of the model can be improved by including higher-order terms as well. The corresponding gSSM-model can be constructed by computing an appropriate Padé-approximant for the forcing term in \eqref{eq:phase_dep}. 

\begin{figure}[h!]
    \centering
    \includegraphics[width = 0.99\textwidth]{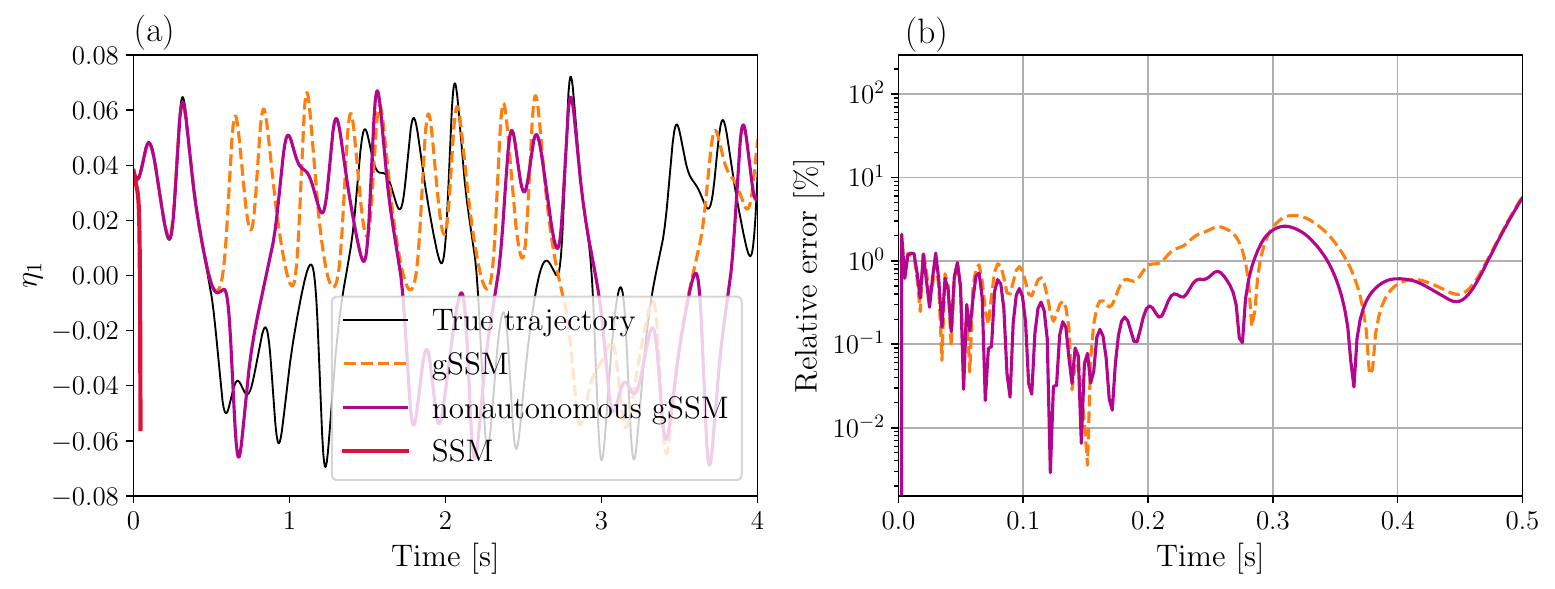}
    \caption{(a): Time series of the reduced coordinate $\eta_1$ on a chaotic trajectory of the full system (black). Higher-order nonautonomous corrections are included in both the SSM-reduced forced model (red) and the gSSM-reduced model (violet). The leading-order approximation of the gSSM-reduced model is also shown orange. (b): Relative error of the reduced-model for short times. }\label{fig:vonkarman_chaotic_nonaut}
\end{figure}

We compare predictions of the leading-order approximation presented in the main text and a [6/6] Padé-approximant computed with $\hat{N}=17$ in Fig. \ref{fig:vonkarman_chaotic_nonaut}. Qualitatively, the same type of chaotic dynamics is observed for the gSSM models, even when higher-order corrections are taken into account. However, for short times, the prediction error decreases even further when the nonautonomous terms are included. We also note that the polynomial SSM-model did not improve, even when we included the phase-dependent terms in \eqref{eq:phase_dep}. Specifically, the model experiences the same finite-time blowup as the leading-order approximation.

\section{Data-driven model of an inverted flag experiment}
In the main text, we have shown that rational function regression is effective in modeling the reduced dynamics on a low-dimensional SSM. Denoting the reduced coordinates by $\boldsymbol{\eta}\in\mathbb{R}^2$, we approximate the reduced dynamics as 

 \begin{equation}
    \label{eq:ratfunc_approx}
    \dot{\boldsymbol{\eta}}(\boldsymbol{\eta}) \approx [N/M](\boldsymbol{\eta}) = \frac{\sum_{|\mathbf{k}|=0}^N \mathbf{a}_{\mathbf{k}}\boldsymbol{\eta}^{\mathbf{k}}}{\sum_{|\mathbf{k}|=0}^M {b}_{\mathbf{k}}\boldsymbol{\eta}^{\mathbf{k}}}.
\end{equation}
In addition, we require that the denominator is non-zero at all points $\boldsymbol{\eta}_i$ in the training set.

\begin{figure}[h!]
    \centering
    \includegraphics[width = 0.99\textwidth]{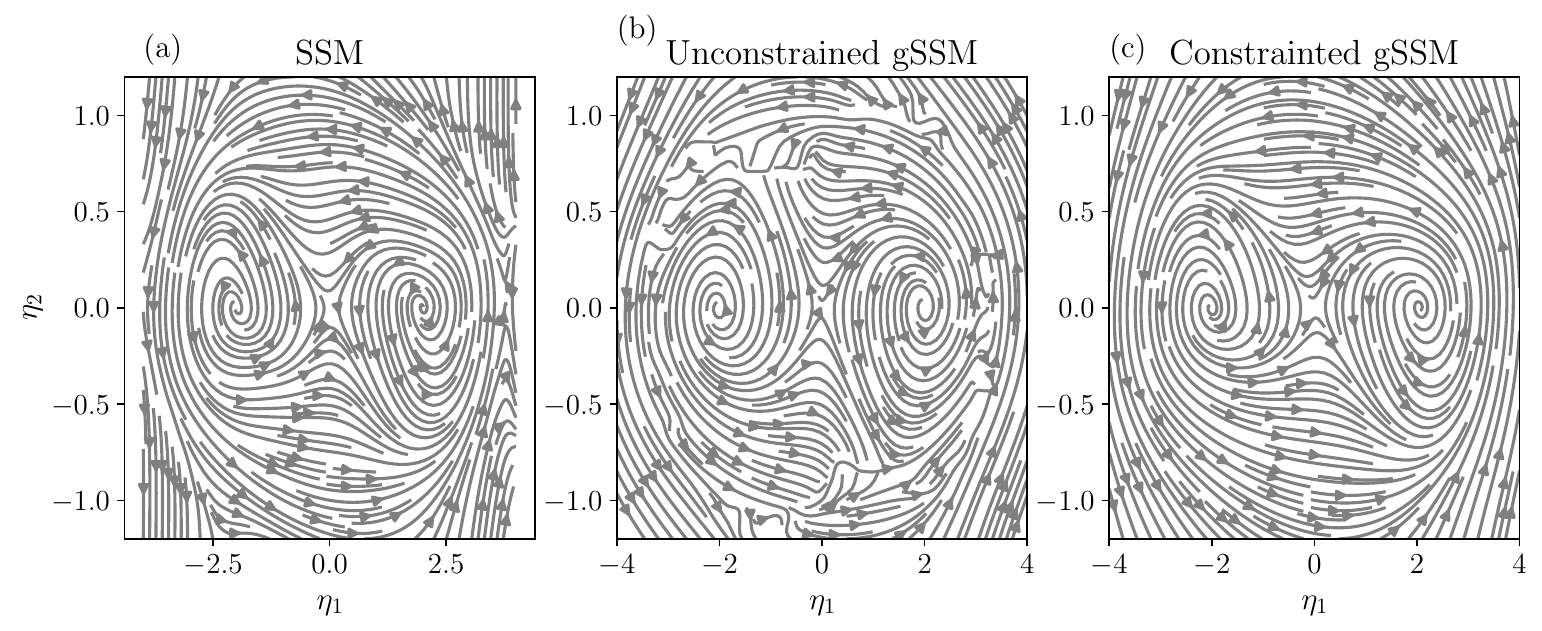}
    \caption{Comparison of data-driven models of the inverted flag experiment. (a): Reduced vector field on the SSM approximated by an order-11 polynomial. (b): The vector field is approximated by a [5/5] rational function without the constraint \eqref{eq:rat_function_2}. (c): Same as (b), but the constraint \eqref{eq:rat_function_2} is enforced.  }\label{fig:comparison_const}
\end{figure}

We then determine the coefficients by minimizing the error 
\begin{equation}
    \label{eq:ratfunc_error}
    \mathcal{E}_r = \sum_{i=1}^K \left| \boldsymbol{\zeta}_i  - \frac{\sum_{|\mathbf{k}|=0}^N \mathbf{a}_{\mathbf{k}}\boldsymbol{\eta}_i^{\mathbf{k}}}{\sum_{|\mathbf{k}|=0}^M {b}_{\mathbf{k}}\boldsymbol{\eta}_i^{\mathbf{k}}} \right|^2,
\end{equation}
such that
\begin{equation}
\label{eq:rat_function_2}
\sum_{|\mathbf{k}|=0}^M {b}_{\mathbf{k}}\boldsymbol{\eta}_i^{\mathbf{k}} \geq \delta \quad \text{ for } i=1,...,K, 
\end{equation}
for some small $\delta>0$.  We point out that without the regularization constraint \eqref{eq:rat_function_2}, the regression can yield spurious singularities, which render the reduced model unusable in practice. In Fig. \ref{fig:comparison_const} we compare the vector fields obtained by polynomial regression (SSM) and rational function regression (gSSM). Due to singularities in the domain of interest, unconstrained rational function regression is unable to recover the correct phase portrait. The polynomial SSM-model and the constrained gSSM-models both capture the dynamical features of the reduced vector field accurately. However, outside the range of the training data bounded by the stable limit cycle, the SSM-model starts to develop large gradients.

\end{document}